# ON A FUNDAMENTAL TASK OF DIOPHANTINE GEOMETRIC FIGURES

Zurab Aghdgomelashvili

(Georgian Technical University, 77, M. Kostava str., Tbilisi, 01785, Georgia)

## Abstract

The goal of the work is to take on and study one of the fundamental tasks studying Diophantine n-gons (the author of the paper considers an integral n-gon is Diophantine as far as determination of combinatorial properties of each of them requires solution of a certain Diophantine equation (equation sets)).

Task*($n; k$): is there a Diophantine n-gon ($n \geq 3$) with any side or diagonal equal to $k$ for each fixed natural number $k$. In case it exists then let us find each such $n$.

It is shown that there is no such Diophantine n-gon ($n > 3$), neither convex nor concave, the length of any side or diagonal of which is equal to 1. It means that for k=1 the above mentioned task is solved.

The studies made it possible to find certain Diophantine rectangles, one of the sides of which is equal to 2 (it is noteworthy that all of them appeared to be inscribed in a circumcircle). The studies showed that there is no such Diophantine rectangle inscribed ina circumcircle, the diagonal length of which is equal to 2. It is shown that for any natural ($k \geq 3$) there is a Diophantine rectangle with the side length equal to $k$.

The fundamental studies showed that for k=2 and n = {3; 4; 5}for convex n-gons (though such Diophantine pentagon has not been found yet. In the author's opinion such pentagon does not exist) and n= {3; 4; 5; 6} for concave n-gons (here as well for n=5 and n=6 no such n-gon has been found yet and for the case of its existence all probable types of such figures are presented).

The paper considers task  *($n; k$) for $k=3$ and shows that in this case 3≤ n ≤7. However, neither such Diophantine pentagon, nor Diophantine hexagon, nor Diophantine heptagon have been found yet. In the author's opinion such Diophantine figures do not exist and in case they do then he presents the probable types of them.

INTRODUCTION

Since ancient times, mathematicians have been interested in the study of geometric figures. It is worthy to mention that great mathematicians, such as: P. Ferma, K.F. Gauss, L. Euler, V. Serpinski, H. Steinhouse et al., were interested in the study of such figures properties. We call such figures as Diophantine due that to determine the properties of each of them they need to be solved by a certain Diophantine equation (system of equations).

By V. Serpinski and H. Steinhouse were stated and solved number of problematic issues related to this topic. Here we will stop on consideration of one of them.

Task: For each $n \in N$ on the plane will be find not located on a line such $n$ points, for that between each two of them the distance will be expressed in natural numbers.

For solving this task, they have shown that for $n \in N$ always exist a Diophantine $n$-gon, the question arises here: exists or not for each fixed $k \in N$ the Diophantine $n$-gon, the distance between any two vertexes of that is equal to $k$, and if so, then find all such $n$.

From the tasks under study on Diophantine geometrical figures, one of the most important is obviously the following:

**Task\* (n;k). Exist or not for any fixed $k$ natural number the Diophantine n-gon ($n \geq 3$) those length of arbitrary side or diagonal is equal to $k$, and if so, then find all such $n$.**

It is shown by us that is not exist such Diophantine $n$-gon, both convex and concave, those lengths of arbitrary side or diagonal are equal to 1. I.e. for $k=1$ the above-mentioned issue has been solved.

BASIC PART

To prove this, let us consider a few tasks.

**Lemma 1. If one of the sides of a Diophantine triangle is equal to 1, then its remaining sides represent equal legs of this triangle.**

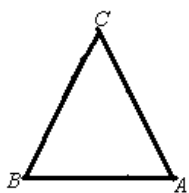
Fig. 1

Task: $\triangle ABC$;
$|BC|, |AC| \in N$;
$|AB| = 1$.
to be proven: $|BC| = |AC|$.

Without limiting the generality let's say $|BC| \leq |AC|$. Due the inequality of triangle $|AB| + |BC| > |AC|$ or $1 + |BC| > |AC|$.

I.e. $\begin{cases} |BC|, |AC| \in N; \\ |BC| \leq |AC| < |BC| + 1. \end{cases} \Rightarrow (|BC| = |AC|).$ Q.E.D..

This remarkable property is one of the cornerstones of the Diophantine and Bidiophantine geometric figure's research apparatus.

**Lemma 2. The length of each side and each diagonal of the convex Diophantine rectangle is greater than 1.**

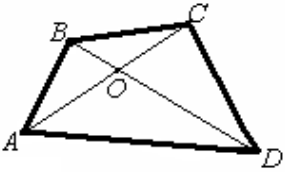

Fig.2

Let's assume the opposite. I.e. let's say that in rectangle *ABCD* the length of each side and diagonal is expressed by the natural number, and also the length of one of the sides, without limiting of generality let's say that $|AB|=1$.

According **to Lemma 1** $|BC|=|AC|$ and $|BD|=|AD|$.

This is impossible, because then *C* and *D* the points should be located on the median $[AB]$. I.e. our assumption is false. Therefore accordingly of this condition, the length of each side of ▱ *ABCD* is greater than 1.

Now let's say that length of the ▱ *ABCD* diagonal is equal to 1, without limitation of generality let's say $|AC|=1$.

According to Lemma 1 $|AB|=|BC|$ and $|CD|=|AD|$. At the same time using the inequality of triangles $\triangle BOC$ and $\triangle AOD$, it is easy to show, that $|BD|+|AC|>|BC|+|AD|$ or $|AB|+|AD|<|BD|+1$. From $\triangle ABD$ we have $|AB|+|AD|>|BD|$. I.e. $\begin{cases} |AB|, |BD|, |AB|\in N; \\ |BD|<|AB|+|AD|<|BD|+1. \end{cases}$ is impossible. Thus the assumption is false, or $|AC|>1$.

**I.e. finally we have that the length of each side and diagonal of each convex Diophantine rectangle is greater than 1. Q.E.D.**

**Lemma 3. If none of the located on the plane four points are not located on one straight line and at the same time distance between each of them is expressed by a natural number, then each of these distances is greater than 1.**

Let's assume the opposite. I.e. let's say that there are such four points on the plane, none of them are located on one straight line, the distance between each of them is expressed by a natural number, and at the same time some of them are equal to 1. Then, according to Lemma 1, the remaining both points must be located on a median of segment with equal to 1 length. Here we would have two cases

1) 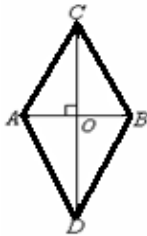 Task: $\triangle ABC$; $|AB|=1$; $|AC|=|BC|\in N$; $|AD|=|DB|\in N$; $|CD|\in N$.

Fig.3

2) 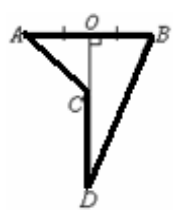 Task: $|AB|=1, |AC|=|BC|=m$; $|AD|=|DB|=n$, $|CD|=l, m,n,l\in N$.

Fig.4

1) This case is considered in **Lemma 2**.

2) From $\triangle AOC$ and $\triangle DOB$ we have:

$$\begin{cases} |AC|^2 = |AO|^2 + |OC|^2; \\ |BD|^2 = |OB|^2 + |OD|^2; \\ |OD| = |OC| + |CD|, |AO| = |OB| = 0,5; \\ |AC| = |CB| = m, |AD| = |DB| = n, |CD| = l, \ m,n,l \in N. \end{cases} \Rightarrow$$

$$\Rightarrow \begin{cases} m^2 = \left(\frac{1}{2}\right)^2 + |OC|^2; \\ n^2 = \left(\frac{1}{2}\right)^2 + (l+|OC|)^2; \\ m,n,l \in N. \end{cases} \Rightarrow \begin{cases} |OC| = \frac{n^2 - m^2 - l^2}{2l}; \\ 4m^2 - 1 = (2|OC|)^2; \\ m,n,l \in N. \end{cases} \Rightarrow \begin{cases} |OC| = \frac{n^2 - m^2 - l^2}{l} = q \in Q_+; \\ 4m^2 - 1 = q^2; \\ m,n,l \in N. \end{cases} \Rightarrow$$

$$\Rightarrow \begin{cases} m,q \in N; \\ 4m^2 - q^2 = 1. \end{cases} \Rightarrow \begin{cases} m,q \in N; \\ (2m-q)(2m+q) = 1. \end{cases} \Rightarrow \begin{cases} m,q \in N; \\ 2m - q = 1; \\ 2m + q = 1. \end{cases}$$

This is impossible. Thus the distance between arbitrary two points from these given points is greater than 1. Q.E.D.

**Theorem 1. The length of each side and each diagonal of convex Diophantine $n$-gon ($n > 3$) is greater than 1.**

Let's assume the opposite. I.e. assume a convex Diophantine $n$-polygon ($n > 3$), the distance between any two vertexes of that is equal to 1. By virtue of **Lemma 1**, the remaining vertexes should be located on the median of the connecting these two vertexes segment. Because this $n$ gon should be convex, so obviously is $n = 4$, but according to Lemma 2 such rectangle does not exist. I.e. our assumption is false, or the length of each side and each diagonal of arbitrary convex Diophantine $n$-gon ($n > 3$) is greater than 1. Q.E.D.

**Theorem 2. If from located on the plane n-point ($n > 3$), none of the three poins is located on one straight line and the distance between each of them is expressed by a natural number, then each of these distances is greater than 1.**

As in the previous theorem, if the distance between arbitrary two vertexes of a Diophantine $n$-polygon is equal to 1, then the its remaining vertexes should be located on the median of connecting these two vertexes segment. If we take into account that none of these vertexes are located on a one straight line, then we will have only the following cases (see Fig. 5).

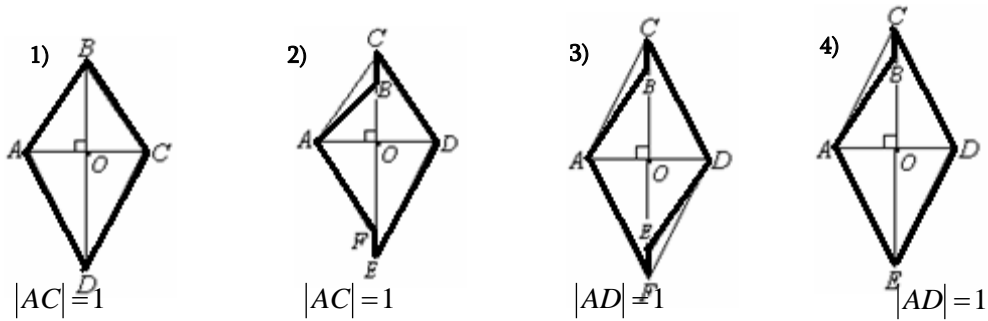

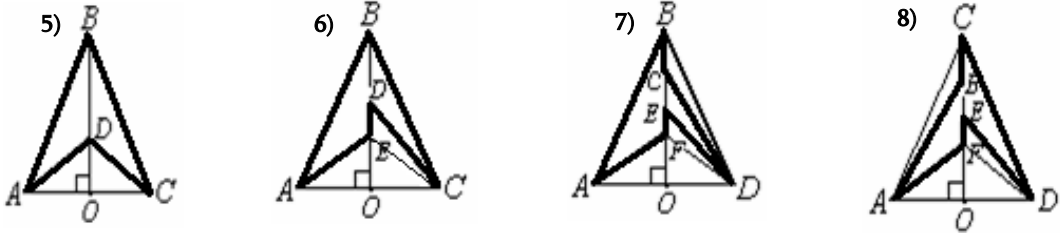

The first fours have no solution according to **Theorem 1.** All fours will be as 1) and because there is not exist such Diophantine rectangle, therefore are not exist 2), 3) and 4) Diophantine rectangles, and the last four have no solution by virtue of **Theorem 2** (Here, too 6), 7) and 8) by fill will be reduced to 5). I.e. our assumption is false or each of the distances given by the task conditions is greater than 1. Thereby fully is proved.

I.e. we have proved the following theorem.

**Task* (n;1) does not exist such a Diophantine n-gon (n≥3), both convex and concave, those lengths of arbitrary sides or diagonals are equal to 1**

Task1. In Diophantine $\triangle ABC$-ძი $|AC|=2$ and $|AB|=a; (a \in N)$ let's find $|BC|$.

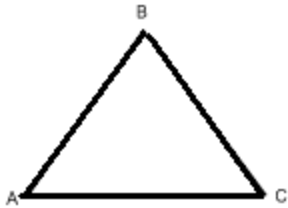

Is given: $\triangle ABC$; $|BC| \in N; |AC| = 2, |AB| = a; a \in N, a \neq 1$.

Must be found. $|BC|$

We have two cases: 1) $\begin{cases} a \in N; a \geq 2; \\ |BC| \leq a. \end{cases}$ 2) $\begin{cases} a \in N; a \geq 2; \\ |BC| \geq a. \end{cases}$

1) From $\triangle ABC$ accordingly of inequality of triangle $a < |BC| + 2$

$$\begin{cases} |BC|; a \in N; a \geq 2; \\ |BC| \leq a < |BC| + 2. \end{cases} \Leftrightarrow \begin{cases} a \in N \\ \left[\begin{array}{l} a = |BC|; \\ a = |BC| + 1. \end{array}\right. \end{cases} \Leftrightarrow \begin{cases} a \in N; a \geq 2; \\ \left[\begin{array}{l} |BC| = a; \\ |BC| = a - 1. \end{array}\right. \end{cases}$$

2) From $\triangle ABC$ accordingly of inequality of triangle $|BC| < a + 2$

$$\begin{cases} |BC|; a \in N; a \geq 2; \\ a \leq |BC| < a + 2. \end{cases} \Leftrightarrow \begin{cases} a \in N; a \geq 2; \\ \left[\begin{array}{l} |BC| = a; \\ |BC| = a + 1. \end{array}\right. \end{cases}$$

By us below are stated all such Diophantine tetragon. length of those side or diagonal would be equal to 2

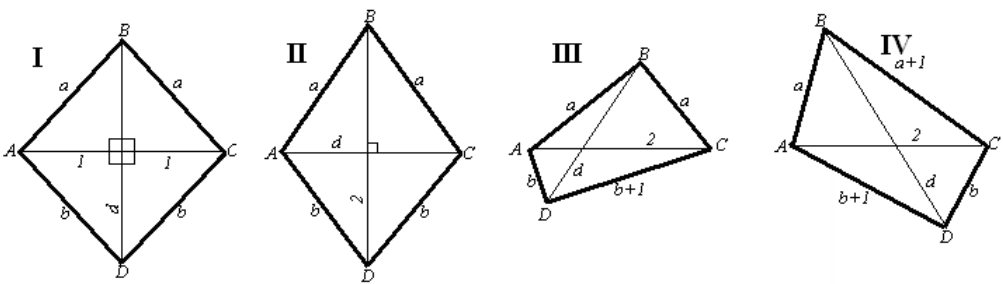

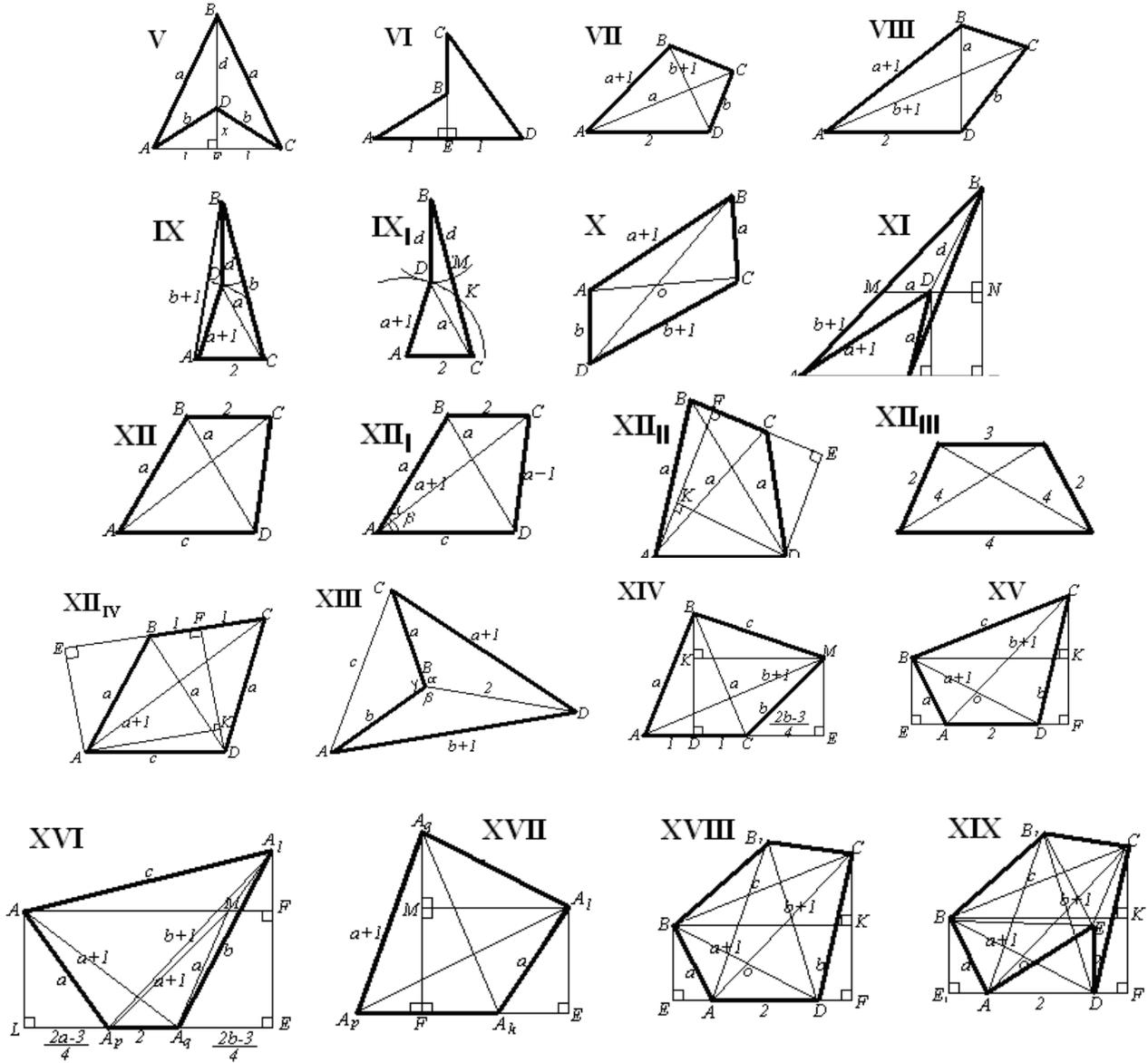

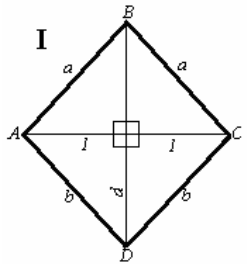

I

Is given: $\begin{cases} a,b,d \in N; \\ |AB|=|BC|=a; \\ |AD|=|CD|=b; \\ |AC|=2; |BD|=d. \end{cases}$

We have already shown in Lemma 3 that such a Diophantine tetragon does not exist.

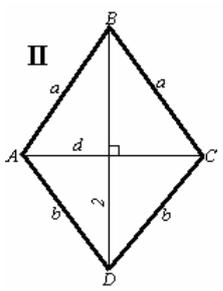

II

Is given: $\begin{cases} a,b,d \in N; a<b; \\ |AB|=|BC|=a; \\ |AD|=|CD|=b; \\ |AC|=d; |BD|=2. \end{cases}$

From $\triangle ABC$ we have $\begin{cases} a,b \in N; \\ a<b; \\ a+2>b. \end{cases} \Rightarrow \begin{cases} a,b \in N; \\ b-2<a<b \end{cases} \Rightarrow (a=b-1)$

It is known that the sum of the diagonals length a convex tetragon is greater than the sum of the lengths of the opposite sides. Thus $d+2 > a+b$ or $d+2 > 2b-1$.

From $\triangle ABC$ $\begin{cases} d,b \in N; \\ d < 2a = 2b-2; \end{cases}$ i.e. $\begin{cases} d,b \in N; \\ 2b-3 < d < 2b-2 \end{cases}$ this is impossible. Thus such Diophantine tetragon is not exist.

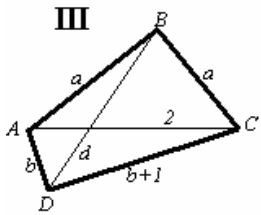

III

Is given: $\begin{cases} a,b,d \in N; |AB| = |BC| = a; \\ |AD| = b; |CD| = b+1 \\ |AC| = 2; |BD| = d. \end{cases}$

From $\triangle ABD$ $d < a+b$.

From the tetragon $ABCD$ $|AC|+|BD| > |AB|+|CD|$, or $d+2 > a+b+1$.

We have: $\begin{cases} d < a+b \\ d+2 > a+b+1 \\ a,b,d \in N \end{cases} \Rightarrow \begin{cases} a+b-1 < d < a+b \\ a,b,d \in N \end{cases}$ this is impossible. Thus such Diophantine tetragon is not exist.

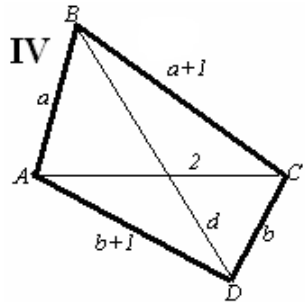

IV

Is given: $\begin{cases} a,b,d \in N; \\ |AB| = a, |BC| = a+1; \\ |AD| = b+1, |CD| = b; \\ |AC| = 2; |BD| = d. \end{cases}$

$\begin{cases} |BD|+|AC| > |BC|+|AD|; \\ |BD| < |AB|+|AD|; \\ |AB|,|BC|,|AD|,|BD| \in N. \end{cases}$ i.e. $\begin{cases} a,b,d \in N; \\ d+2 > a+1+b+1 \\ d < a+b+1 \end{cases} \Rightarrow \begin{cases} a,b,d \in N; \\ a+b < d < a+b+1 \end{cases}$

This is impossible. Thus such Diophantine tetragon is not exist.

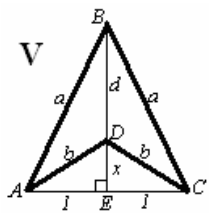

V

Is given: $\begin{cases} a,b,d \in N; \\ |AB| = |BC| = a; \\ |AD| = |CD| = b; \\ |BD| = d; |AC| = 2 \end{cases}$

From the rectangle $\triangle AED$ and rectangle $\triangle AEB$ we have:

$\begin{cases} a,b,d \in N \\ 1^2 + x^2 = b^2 \\ 1^2 + (d+x)^2 = a^2 \end{cases} \Rightarrow \begin{cases} a,b,d \in N \\ 1+x^2 = b^2 \\ 1+x^2+d^2+2dx = a^2 \end{cases} \Rightarrow \begin{cases} a,b,d \in N \\ x = \dfrac{a^2-b^2-d^2}{2d} \in Q_+ \\ 1+x^2 = b^2 \end{cases} \Rightarrow \begin{cases} b,x \in N \\ (b-x)(b+x) = 1 \end{cases} \Rightarrow$

$\Rightarrow \begin{cases} b,x \in N \\ b-x = 1 \\ b+x = 1 \end{cases}$ This is impossible. Thus such Diophantine tetragon is not exist.

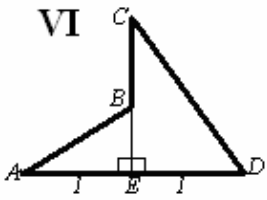

Is given: $\begin{cases} |AB| = |BD| = a; |AC| = |CD| = b; \\ |BC| = d; |AD| = 2. \\ a,b,d \in N \end{cases}$

VI

Proceeding from the task conditions the $B$ and $C$ points are located on the middle perpendicular of $[AD]$, thus $|AE| = |ED| = \dfrac{|AD|}{2} = 1$.

From the rectangle $\triangle AEB$ and rectangle $\triangle CED$ we have:

$\begin{cases} a^2 = |BE|^2 + 1 \\ b^2 = (d + |BE|)^2 + 1; \\ a,b,d \in N. \end{cases} \Rightarrow \begin{cases} a^2 = |BE|^2 + 1 \\ |BE| = \dfrac{b^2 - d^2 - a^2}{2d} \in Q_+ \\ a,b,d \in N. \end{cases} \Rightarrow \begin{cases} a^2 - |BE|^2 = 1 \\ a, |BE| \in N \end{cases} \Rightarrow$

$\Rightarrow \begin{cases} (a - |BE|)(a + |BE|) = 1; \\ a, |BE| \in N. \end{cases} \Rightarrow \begin{cases} a - |BE| = 1; \\ a + |BE| = 1; \\ a, |BE| \in N. \end{cases}$ This is impossible. Thus such Diophantine tetragon is not exist.

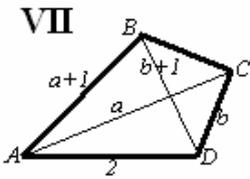

VII

Is given: $\begin{cases} a,b,d,c \in N; \\ |AB| = a+1; |AC| = a; |BC| = c; |CD| = b; \\ |BD| = b+1; |AD| = 2. \end{cases}$

Such Diophantine tetragon is not exist, because the sum of the lengths of diagonal of the convex tetragon will be greater than the sum of the lengths of the opposite sides of this tetragon, i.e. it should be $|BD| + |AC| > |AB| + |CD|$, but proceeding from this we would have $a + (b+1) > (a+1) + b$, that is impossible.

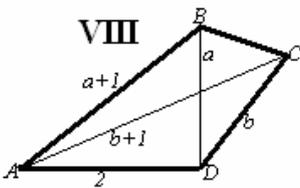

VIII

Is given: $\begin{cases} a,b,d,c \in N; \\ |AB| = a+1; |AC| = b+1; |BC| = c; |CD| = b; \\ |BD| = a; \ |AD| = 2. \end{cases}$

In VIII, similarly to VII would be $|BD| + |AC| > |AB| + |CD|$, but proceeding from this we would have $a + (b+1) > (a+1) + b$. Thus such Diophantine tetragon is not exist.

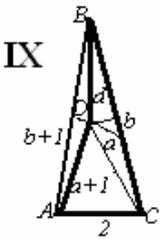 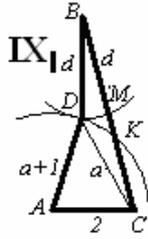

IX

Is given: $\begin{cases} a,b,d,c \in N; \\ |AB| = b+1; |AC| = 2; |BC| = b; |CD| = a; \\ |BD| = d; \ |AD| = a+1. \end{cases}$

Let's circumscribe from the point A by radius $[AD]$ and point B by radius $[BD]$ the arcs and say they intersect $[BC]$ at points K and M, respectively.

From $\triangle AKC$ we have $(2 + |CK| > a+1) \Rightarrow (|CK| > a - 1)$.

$$\begin{cases} |CK|+|BM|<b; \\ |CK|>a-1; |BM|=d. \end{cases} \Rightarrow a-1+d<b \Rightarrow d<b-a+1.$$

$\Delta BDC$-დან $d+a>b \Rightarrow d>b-a$.

i.e. $\begin{cases} d>b-a; \\ d<b-a+1 \\ a,b,d \in N \end{cases} \Rightarrow \begin{cases} b-a<d<b-a+1; \\ a,b,d \in N \end{cases}$ This is impossible. Thus such Diophantine tetragon is not exist.

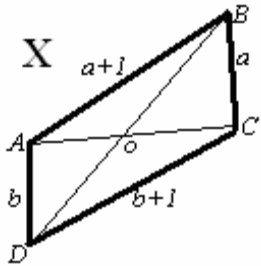

Is given: $\begin{cases} a,b,|BD| \in N; a>3; b>3; \\ |AB|=a+1; |BC|=a; |AD|=b; |CD|=b+1; \\ \left[\begin{array}{l}|AC|=2;\\|AC|=3.\end{array}\right. \end{cases}$

$\begin{cases} a,b \in N \\ a>3 \\ b>3 \\ \left[\begin{array}{l}|AC|=2\\|AC|=3\end{array}\right. \end{cases} \Rightarrow \begin{cases} |AC|^2+a^2 \le (a+1)^2 \\ |AC|^2+b^2 \le (b+1)^2 \end{cases} \Rightarrow \begin{cases} A\widehat{C}B \in \left[\dfrac{\pi}{2};\pi\right) \\ D\widehat{A}C \in \left[\dfrac{\pi}{2};\pi\right) \end{cases} \Rightarrow \begin{cases} a<|BD|<a+1 \\ b<|DO|<b+1 \end{cases} \Rightarrow$

$\Rightarrow a+b<|BO|+|OD|<a+b+2 \Rightarrow |BD| \in (a+b; a+b+2)$.

From $\Delta BCD$ we have $|BD|<a+b+1$. i.e. $\begin{cases} |BD| \in (a+b; a+b+1); \\ |BD| \in N \end{cases} \Rightarrow |BD| \notin N$.

Thus such Diophantine tetragon is not exist.

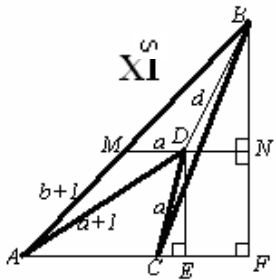

$\Delta MBN \sim \Delta ABF \Rightarrow \dfrac{|MB|}{|AB|}=\dfrac{|BN|}{|BF|} \Rightarrow$

$\Rightarrow \dfrac{|MB|}{b+1}=\dfrac{\sqrt{12b^2+12b-9}-\sqrt{12a^2+12a-9}}{\sqrt{12b^2+12b-9}}=1-\sqrt{\dfrac{2a-1}{2b-1} \cdot \dfrac{2a+3}{2b+3}}<$

$<1-\dfrac{2a-1}{2b-1}=\dfrac{b-a}{b-\dfrac{1}{2}}$.

We have $|MB|<(b-a)+1.5\left(\dfrac{b-a}{b-\dfrac{1}{2}}\right)<(b-a)+1.5$.

$\begin{cases} |BD|<|MB|<(b-a)+1.5 \\ |BD| \in N \end{cases} \Rightarrow |BD| \le (b-a)+1$

From $\triangle ADB$ we have: $(a+1)+|BD| > b+1 \Rightarrow |BD| > b-a$.

i.e. $\begin{cases} b-a < |BD| \leq b-a+1; \\ a,b,|BD| \in N \end{cases} \Rightarrow |BD| = b-a+1$.

From the rectangle $\triangle DNB$ we have $|DN|^2 + |BN|^2 = |BD|^2$. i.e.

$$\left(\frac{b-a}{2}\right)^2 + \left(\frac{\sqrt{12b^2+12b-9}-\sqrt{12a^2+12a-9}}{4}\right)^2 = ((b-a)+1)^2.$$ From that we will obtain:

$$24a^2b - 24ab^2 - 62ab + 37a^2 + 13b^2 + 28b - 40a + 13 = 0 \quad (*)$$

Proceeding from this

$$37a^2 + 13b^2 + 13 \equiv 0 \pmod 2 \Rightarrow a^2 + b^2 + 1 \equiv 0 \pmod 2 \Rightarrow \begin{cases} a=2p \\ b=2k-1 \\ p,k \in N \end{cases} \cup \begin{cases} a=2p-1 \\ b=2k \\ p,k \in N \end{cases}$$

If $b=2k$ then from (*) we have: $37a^2 + 13 = 4\cdot 37p(p-1) + 37 + 13 \equiv 0 \pmod 4$, that is impossible.

if $a=2p$, then from (*) we have $13b^2 + 13 = 13(4p(p-1)+2) \equiv 0 \pmod 4$, that also is impossible. I.e. such Diophantine rectangle is not exist.

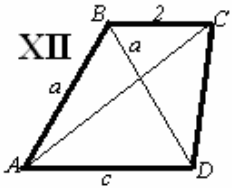

**XII** Is given: $\begin{cases} |AB|=|BD|=a, \ |AD|=c, \ |BC|=2; \\ |AC|, |CD|, a, c \in N; |CD| < |BD|. \end{cases}$

From the conditions of task we have: $c<2a$;

$|CD| = a-1;$ $\begin{bmatrix} |AC| = a+1; & (*) \\ |AC| = a & (**) \end{bmatrix}$

Let's consider each of them separately:

(*) 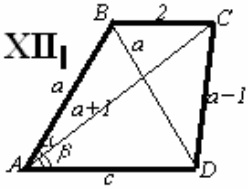

From $\triangle ABC$ and $\triangle ACD$ due the cosine theory we have:

$$\begin{cases} 2^2 = a^2 + (a+1)^2 - 2a(a+1)\cos\alpha \\ (a-1)^2 = c^2 + (a+1)^2 - 2c(a+1)\cos\beta \end{cases} \Rightarrow$$

$$\begin{cases} \cos\alpha = \dfrac{2a^2+2a-3}{2}; \ \sin\alpha = \dfrac{\sqrt{12a^2+12a-9}}{2a(a+1)}; \\ \cos\beta = \dfrac{c^2+4a}{2c(a+1)}; \ \sin\beta = \dfrac{\sqrt{(c^2-4)(4a^2-c^2)}}{2c(a+1)} \end{cases}$$

From $\triangle ABD$

$$\frac{c}{2a} = \cos(\alpha+\beta) = \frac{(c^2+4a)(2a^2+2a-3)}{4ac(a+1)^2} - \frac{\sqrt{(c^2-4)(4a^2-c^2)(12a^2+12a-9)}}{4ac(a+1)^2} \Rightarrow$$

$$\Rightarrow \sqrt{(c^2-4)(4a^2-c^2)(12a^2+12a-9)} = -(2a+5)c^2 + (8a^3+8a^2-12a) \Rightarrow$$

$$\Rightarrow 16(a+1)^2 c^4 - 4(a+1)^2(20a^2-9)c^2 + 64a^4(a+1)^2 = 0 \Rightarrow 4c^4 - (20a^2-9)c^2 + 16a^4 = 0.$$

Proceeding from this $c = 4p, (p \in N)$. Due the introduction we will obtain:

$$\left(4(4p)^4 - (20a^2-9)(4p)^2 + 16a^4 = 0\right) \Leftrightarrow \left(64p^4 - (20a^2-9)p^2 + a^4 = 0\right)$$

From the last equation we have: $\begin{cases} a^4 = 0 \pmod{p^2} \\ a, p \in N \end{cases} \Rightarrow \begin{cases} a^2 = pt \\ a, p, t \in N \end{cases}$ Due the introduction of

obtained in the last equation and it's simplification we will obtain:

$$\begin{cases} t^2 - 20pt + (64p^2+9) = 0 \\ p, t \in N \end{cases} \Leftrightarrow \begin{cases} \left[\begin{array}{l} t = 10p - 3\sqrt{4p^2-1}; \\ t = 10p + 3\sqrt{4p^2-1}. \end{array}\right. \Rightarrow t \in \varnothing. \\ p, t \in N \end{cases}$$

**Remark**: $\begin{cases} 4p^2 - 1 = n^2; \\ p, n \in N \end{cases} \Leftrightarrow \begin{cases} (2p-n)(2p+n) = 1; \\ p, n \in N \end{cases} \Leftrightarrow \begin{cases} 2p - n = 1; \\ 2p + n = 1; \Rightarrow (n \in \varnothing). \\ p, n \in N. \end{cases}$

i. e. such Diophantine tetragon is not exist.

(**)

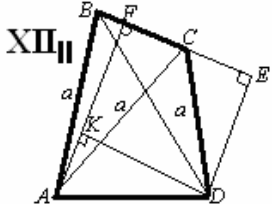

Is given: $\begin{cases} |AB| = |BD| = |AC| = a; \\ |CD| = a-1; |BC| = 2; |AD| = c; \\ a, c \in N. \end{cases}$

From the rectangle $\Delta BED$ and rectangle $\Delta CED$ we have:

$$a^2 - (2+|CE|)^2 = (a-1)^2 - |CE|^2 \Rightarrow |CE| = \frac{2a-5}{4};$$

$$|DE| = \sqrt{(a-1)^2 - \left(\frac{2a-5}{4}\right)^2} = \frac{\sqrt{12a^2-12a-9}}{4}$$

From the rectangle $\Delta AKD$ due the Pythagorean theorem:

$$c^2 = \left(\frac{2a-1}{4}\right)^2 + \left(\sqrt{a^2-1} - \frac{\sqrt{12a^2-12a-9}}{4}\right)^2 \Rightarrow 4c^4 - 4(4a^2-2a-3)c^2 + (4a^4-4a^3+a^2) = 0$$

where $a, c \in N$, $a > 2$, $c > 1$.

From the obtained equation leads that: $a=2p; p \in N$. By introduction we will obtain:

$$\begin{cases} c^4 - (16p^2-4p-3)c^2 + p^2(4p-1)^2 = 0 \\ p, c \in N, p, c > 1 \end{cases} \quad (1)$$

From (1) for $p=2$ we will have:

$$\begin{cases} c^4 - 53c^2 + 192 = 0 \\ c \in N, c > 1 \end{cases} \Leftrightarrow \left[\begin{array}{l} c = 2 \\ c = 7. \end{array}\right.$$

For c=7 the $\triangle ACD$ is degenerated, because $|AC|+|CD|=4+3=7=|AD|$.

For c=2 we will obtain the yet known to us Diophantine trapezoid. The arbitrary other solution of (1) currently is unknown. In our other natural solution is not exist.

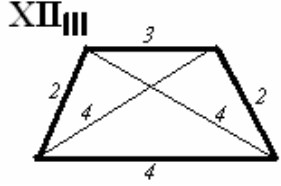
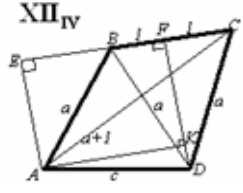

Is given: $\begin{cases} |AB|=|BD|=|CD|=a; \\ |AC|=a+1; |BC|=2; |AD|=c; \\ a,c \in N. \end{cases}$

From the rectangle $\triangle AEC$ and rectangle $\triangle AEB$ we have:

$$(a+1)^2 - (|EB|+2)^2 = a^2 - |EB|^2 \Rightarrow |EB| = \frac{2a-3}{4}; |AE| = \sqrt{a^2 - \left(\frac{2a-3}{4}\right)^2} = \frac{\sqrt{12a^2+12a-9}}{4}.$$

From the rectangle $\triangle AKD$ due the Pythagorean Theorem:

$$c^2 = \left(\frac{2a+1}{4}\right)^2 + \left(\sqrt{a^2-1} - \frac{\sqrt{12a^2+12a-9}}{4}\right)^2,$$

from that we will obtain: $4c^4 - 4(4a^2+2a-3)c^2 + a^2(2a+1)^2 = 0$. Similarly to XII$_{III}$ will be obtained:

$$\begin{cases} c^4 - (16p^2+4p-3)c^2 + p^2(4p+1)^2 = 0 \\ p,c \in N, p,c > 1 \\ a = 2p \end{cases} \quad (2)$$

Currently is not found even one solution of (2). In our opinion it will not have the solution on natural numbers.

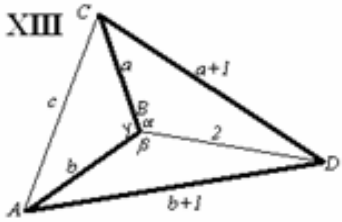

Is given: $\begin{cases} \alpha + \beta > 180°; \\ |BC|=a; |CD|=a+1; |AD|=b+1; \\ |AB|=b; |BD|=2; |AC|=c; \\ a,b,c \in N. \end{cases}$

From $\triangle BCD$; $\triangle ABD$ and $\triangle ABC$ due the cosine's theorem we have:

$\begin{cases} (a+1)^2 = a^2 + 2^2 - 4a \cdot \cos\alpha; \\ (b+1)^2 = b^2 + 2^2 - 4b\cos\beta; \\ c^2 = a^2 + b^2 - 2ab\cos\gamma; \\ \cos\gamma = \cos(360° - (\alpha+\beta)) = \cos(\alpha+\beta). \end{cases} \Rightarrow \begin{cases} \cos\alpha = \frac{3-2a}{4a}; \sin\alpha = \frac{\sqrt{12a^2+12a-9}}{4a} \\ \cos\beta = \frac{3-2b}{4b}; \sin\beta = \frac{\sqrt{12b^2+12b-9}}{4b} \\ \cos\gamma = \cos(\alpha+\beta) = \frac{a^2+b^2+c^2}{2ab}. \end{cases} \Rightarrow$

$$\Rightarrow \frac{a^2+b^2-c^2}{2ab} = \cos(\alpha+\beta) = \frac{(3-2a)(3-2b) - \sqrt{(12a^2+12a-9)(12b^2+12b-9)}}{16ab} \Rightarrow$$

$$\Rightarrow \sqrt{(12a^2+12a-9)(12b^2+12b-9)} = (3-2a)(3-2b) - 8(a^2+b^2-c^2) \Leftrightarrow$$

$$\Leftrightarrow 4c^4 + c^2(9 - 6b - 6a + 4ab - 8a^2 - 8b^2) + 2(a-b)^2(2(a^2 + ab + b^2) + 3(a+b)) = 0 \quad (1)$$

If we in (1) introduce $a = b - b$, we will obtain:

$$4c^4 + c^2(9 - 6a - 6a + 4a^2 - 8a^2 - 8a^2) + 2(a-a)^2(2(a^2 + a^2 + a^2) + 3(a+a)) = 0 \Leftrightarrow$$

$\Leftrightarrow 4c^2 = 12a^2 + 12a - 9$, Which has no solution in natural numbers, because otherwise its left side would be flat and the right side would be odd. As for (1), at present we do not know at least one in the natural solutions, because, in contrary, it's left side will be even and the rigjt side – the odd. As for (1), Currently for us is unknown even one from triplet natural solutions.

**XIV**

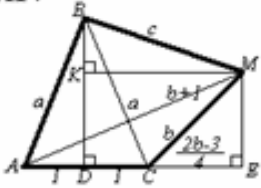

Is given: $\begin{cases} |AB| = |BC| = a; |MC| = b; \\ |AM| = b+1; |BM| = c; |AC| = 2; \\ a, b, c \in N. \end{cases}$

$$|MK| = 1 + \frac{2b-3}{4} = \frac{2b+1}{4};$$

$$|BK| = \sqrt{c^2 - \left(\frac{2b+1}{4}\right)^2}; |KD| = |ME| = \frac{\sqrt{12b^2 + 12b - 9}}{4};$$

$$\begin{cases} |BD| = \sqrt{4c^2 - 1}; \\ |BD| = \sqrt{c^2 - \left(\frac{2b+1}{4}\right)^2} + \frac{\sqrt{12b^2 + 12b - 9}}{4} \end{cases} \Rightarrow 4\sqrt{4c^2 - 1} = \sqrt{16c^2 - (2b+1)^2} + \sqrt{12b^2 + 12b - 9} \Leftrightarrow$$

$$\Leftrightarrow 9c^4 - 6c^2(b^2 + b) + (b^2 + b)^2 = 0 \Leftrightarrow$$

$$\Leftrightarrow \left(\frac{3c^2}{b^2 + b} - 1\right)^2 = 0 \Leftrightarrow \frac{3c^2}{b^2 + b} - 1 = 0 \Leftrightarrow b^2 + b = 3c^2 \Leftrightarrow b(b+1) = 3c^2 \Leftrightarrow (2b+1)^2 - 12c^2 = 1. \quad (1)$$

$$\left(a_1^2 - 12b_1^2\right)\left(c_1^2 - 12d_1^2\right) = \left(a_1 c_1 - 12 b_1 d_1\right)^2 - 12\left(a_1 d_1 - b_1 c_1\right)^2. \quad (2)$$

(1) represents the Pelli equation. By it's solution we obtain the natural solutions of (1). Amongst them first three pairs are pairs:

$$\left[\begin{cases} 2b+1 = 7; \\ c = 2. \end{cases}; \begin{cases} 2b+1 = 97; \\ c = 28. \end{cases}; \begin{cases} 2b+1 = 1351; \\ c = 390. \end{cases}\right] \Leftrightarrow \left[\begin{cases} b = 3; \\ c = 2. \end{cases}; \begin{cases} b = 48; \\ c = 28. \end{cases}; \begin{cases} b = 675; \\ c = 390. \end{cases}\right]$$

**XV**

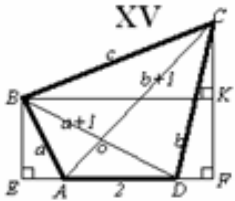

Is given: $\begin{cases} AB = a; |BC| = c; |CD| = b; |AC| = b+1; \\ |BD| = a+1; a, b, c \in N. \end{cases}$

$\triangle ACD$ -ො ხ $\cos A\widehat{D}C = \dfrac{3-2b}{4b} < 0 \; (b \in N; b \neq 1).$

Accordingly of one of the Kumer theorem(**theorem: if the lengths of the sides and diagonals of a convex tetragon are expressed by rational numbers, then the diagonals at the intersection point are divided by rational length sections**) $|AO|, |OD| \in Q_+$.

From $\triangle AOD$ due the law of sines we will have:

$$\frac{|AO|}{|OD|} = \frac{\sin O\widehat{D}A}{\sin O\widehat{A}D} = \frac{\sin B\widehat{D}A}{\sin C\widehat{A}D} = \frac{\sqrt{1-\cos^2 B\widehat{D}A}}{\sqrt{1-\cos^2 C\widehat{A}D}} = \frac{\sqrt{1-\left(\frac{2a+5}{4(a+1)}\right)^2}}{\sqrt{1-\left(\frac{2b+5}{4(b+1)}\right)^2}} = \frac{(b+1)\sqrt{4a^2+4a-3}}{(a+1)\sqrt{4b^2+4b-3}}.$$

Proceeding from this $\dfrac{\sqrt{4a^2+4a-3}}{\sqrt{4b^2+4b-3}} = t \in Q_+$ or $4a^2+4a-3 = t^2(4b^2+4b-3)$ .i.e.this solution of takes is reduced on the solution of following task.

Let's solve in natural $a$, $b$, $m$ и $n$ numbers the equation

$$n^2(4a^2+4a-3) = m^2(4b^2+4b-3).$$

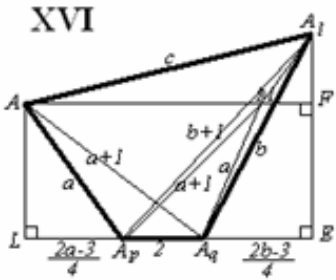

XVI

$$|AM| = \frac{2a-3}{2} + 2 = \frac{2a+1}{2}$$

From the tetragon $A_p A_q A_l M$

$$|AM| + |A_p A_l| > |AA_l| + A_p M \Rightarrow$$

$$\Rightarrow \frac{2a+1}{2} + b + 1 > a + 1 + c \Rightarrow c < b + \frac{1}{2} \Rightarrow c \leq b.$$

$b > a \Rightarrow \dfrac{2a-1}{2b-1} < \dfrac{2a+3}{2b+3}$ From the right triangle $\Delta AFA_l$ due the Pythagorean Theorem:

$$c^2 = \left(\frac{a+b+1}{2}\right)^2 + \frac{1}{16}\left(\sqrt{12b^2+12b-9} - \sqrt{12a^2+12a-9}\right)^2 =$$

$$= \left(\frac{a+b+1}{2}\right)^2 + \frac{12b^2+12b-9}{16}\left(1 - \sqrt{\frac{2a-1}{2b-1} \cdot \frac{a+3}{2b+3}}\right)^2 = \left(\frac{a+b+1}{2}\right)^2 + \frac{3(2b-1)(2b+3)}{16} \times$$

$$\times \left(1 - \sqrt{\frac{2a-1}{2b-1} \cdot \frac{a+3}{2b+3}}\right)^2 \quad (*).\text{ From } (*)\text{ we have:}$$

$\boxed{1}$ $\quad c^2 < \left(\dfrac{a+b+1}{2}\right)^2 + \dfrac{3(2b-1)(2b+3)}{16} \cdot \dfrac{4(b-a)^2}{(2b-1)^2} = \left(\dfrac{a+b+1}{2}\right)^2 + \dfrac{3}{4}\left(\dfrac{2b+3}{2b-1}\right) \cdot (b-a)^2;$

$\boxed{2}$ $\quad c^2 > \left(\dfrac{a+b+1}{2}\right)^2 + \dfrac{3(2b-1)(2b+3)}{16} \cdot \dfrac{4(b-a)^2}{(2b-3)^2} = \left(\dfrac{a+b+1}{2}\right)^2 + \dfrac{3}{4}\left(\dfrac{2b-1}{2b+3}\right) \cdot (b-a)^2.$

I.e. if exists such Diophantine tetragon, then $c$ will satisfy the conditions $\boxed{1}$ and $\boxed{2}$.

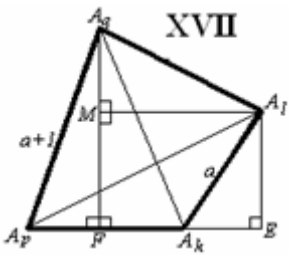

XVII

Is given: $\begin{cases} A_p A_q A_l A_k; \\ |A_p A_q| = |A_q A_k| = |A_p A_l| = a+1 \\ |A_k A_l| = a; |A_p A_k| = 2; \\ |A_q A_l| = c; \\ a, c \in N. \end{cases}$

From the rectangles $\Delta A_p E A_l$ and $\Delta A_k E A_l$ we have:

$$(a+1)^2 - (|A_k E| + 2)^2 = a^2 - |A_k E|^2 \Rightarrow a^2 + 2a + 1 - |A_k E|^2 - 4|A_k E| - 4 = a^2 - |A_k E|^2 \Rightarrow$$

$$\Rightarrow |A_k E| = \frac{2a-3}{4}.$$

$$|A_q F| = \sqrt{(a+1)^2 - 1} = \sqrt{a^2 + 2a}; \quad |EF| = \frac{2a-3}{4} + 1 = \frac{2a+1}{4}.$$

$$|A_l E| = \sqrt{a^2 - \left(\frac{2a-3}{4}\right)^2} = \frac{\sqrt{12a^2 + 12a - 9}}{4}.$$

From $\Delta A_q M A_l$ due the Pythagorean Theorem we have:

$$c^2 = \left(\frac{2a+1}{4}\right)^2 + \left(\sqrt{a^2 + 2a} - \frac{\sqrt{12a^2 + 12a - 9}}{4}\right)^2 = \frac{4a^2 + 6a - 1 - \sqrt{12a^4 + 36a^3 + 15a^2 - 18a}}{2}.$$

Due the simple transformations we will obtain the equation:

$$c^2 = \frac{4a^2 + 6a - 1 - \sqrt{(a^2 + 2a)(12a^2 + 12a - 9)}}{2}. \qquad (1)$$

$\boxed{1}$ $a \geq 3 \Leftrightarrow 12a^2 + 12a - 9 \geq 9(a^2 + 2a) \Leftrightarrow \sqrt{(a^2 + 2a)(12a^2 + 12a - 9)} \geq \sqrt{9(a^2 + 2a)^2} \Leftrightarrow$

$$\Leftrightarrow c^2 \leq \frac{4a^2 + 6a - 1 - 3(a^2 + 2a)}{2} = \frac{a^2 - 1}{2} \Leftrightarrow c \leq \sqrt{\frac{a^2 - 1}{2}}$$

$\boxed{2}$ $\quad a^2 + 50a + 36 > 0 \Leftrightarrow 48a^2 + 48a - 36 < 49a^2 + 98a \Leftrightarrow$

$$\Leftrightarrow 12a^2 + 12a - 9 < \frac{49}{4}(a^2 + 2a) \Leftrightarrow (a^2 + 2a)(12a^2 + 12a - 9) < 3.5^2 (a^2 + 2a)^2 \Leftrightarrow$$

$$\Leftrightarrow \sqrt{(a^2 + 2a)(12a^2 + 12a - 9)} < 3.5(a^2 + 2a)^2.$$

I.e. $c^2 > \dfrac{4a^2 + 6a - 1 - 3.5(a^2 + 2a)}{2} = \dfrac{0.5a^2 - a - 1}{2} = \left(\dfrac{a-1}{2}\right)^2 - \dfrac{3}{4} \Rightarrow c \geq \dfrac{a-1}{2}.$

Finally we have that if exists such Diophantine tetragon, then $c \in \left[\dfrac{a-1}{2}; \sqrt{\dfrac{a^2-1}{2}}\right)$.

Let's now consider the issue of inscribed in a circle and circumscribed on a circle for such Diophantine tetragons, the length of one of those sides is equal to 2.

As we have seen, in some cases on the XIV type Diophantine tetragon will be circumscribed the circle. Now let's study whether or no inscribe in this kind of Diophantine tetragonthe circle?

**Task2**.

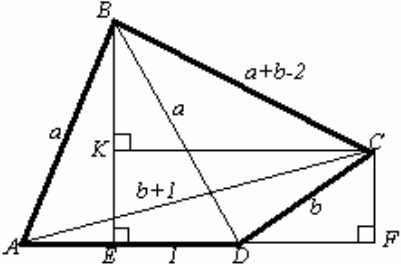

Is given: In ▱ABCD will be inscribed the circle

$$\begin{cases} |AB|=|BD|=a;\ |CD|=b;\\ |AC|=b+1;\\ |AD|=2\\ a,b,|BC|\in N. \end{cases}$$

(In ▱ABCD will be inscribed the circle) $\Rightarrow |BC|=|AB|+|CD|-|AD|=a+b-2$.

$|BD|=\sqrt{|BE|^2+|ED|^2}\Rightarrow |BE|=\sqrt{a^2-1}$

$|BE|=|BK|+|KE|=\sqrt{|BC|^2-|KC|^2}+|CF|=\sqrt{(a+b-2)^2-\left(\dfrac{2b+1}{4}\right)^2}+\dfrac{\sqrt{12b^2+12b-9}}{4}$.

I.e. $\sqrt{a^2-1}=\sqrt{(a+b-2)^2-\left(\dfrac{2b+1}{4}\right)^2}+\dfrac{\sqrt{12b^2+12b-9}}{4}\Leftrightarrow$

$\Leftrightarrow \sqrt{16(a+b-2)^2-(2b+1)^2}=\sqrt{16a^2-16}-\sqrt{12b^2+12b-9}\Leftrightarrow$

$\Leftrightarrow \begin{cases}\sqrt{16a^2-16}-\sqrt{12b^2+12b-9}>0\\ 16a^2+12b^2+32ab-64a-68b+63=16a^2-16+12b^2+12b-9-8\sqrt{(a^2-1)(12b^2+12b-9)}\end{cases}\Leftrightarrow$

$\Leftrightarrow \begin{cases}16a^2>12b^2+12b+7\\ -4ab+8a+10b-11=\sqrt{(a^2-1)(12b^2+12b-9)}.\end{cases}$ (1)

$\begin{cases}a>2;\\ b>2.\end{cases}\Rightarrow -4ab+8a+108-11\geq \sqrt{(2^2-1)(12\cdot 2^2+12\cdot 2-9)}=\sqrt{189}>13\Rightarrow 2a(b-2)<5b-12$ (*)

$(b>2)\Rightarrow a<2.5-\dfrac{1}{b-2}<2.5.$ i.e. $a=2$.

Due the introduction of $a=2$ in (1) we will obtain: $\begin{cases}-8b+16+10b-11=3\sqrt{4b^2+4b-3}\\ b\in N.\end{cases}\Leftrightarrow$

$\Leftrightarrow \begin{cases}2b+5=3\sqrt{4b^2+4b-3}\\ b\in N\end{cases}\Leftrightarrow \begin{cases}8b^2+4b-13=0\\ b\in N\end{cases}\Leftrightarrow b\in\varnothing$

Due the introduction of $b=2$ in (1) we will obtain: $\begin{cases}-8a+8a+20-11=3\sqrt{7(a^2-1)}\\ a\in N;\ a\geq 2.\end{cases}\Leftrightarrow$

$\Leftrightarrow \begin{cases}9=7(a^2-1)\\ a\in N;\ a\geq 2.\end{cases}\Rightarrow a\in\varnothing$. I.e. such Diophantine tetragon is not exist.

**Let's now consider the issue of inscribed in a circle and circumscribed on a circle for XIV type Diophantine tetragons.**

Task3.

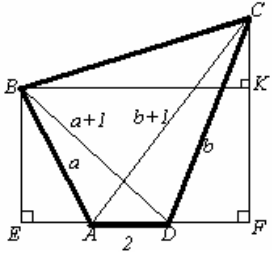

Is given: *ABCD* is inscribed in the circle;

$$\begin{cases} |AB| = a; |BD| = a+1; |CD| = b; \\ |AC| = b+1; |AD| = 2; |BC| = c; \\ a, b, c \in N \end{cases}$$

(□ *ABCD* is inscribed in the circle)

$$\Rightarrow |BD|\cdot|AC| = |AB|\cdot|CD| + |BC|\cdot|AD| \Rightarrow (a+1)(b+1) = ab + 2c \Leftrightarrow c = \frac{a+b+1}{2};$$

From the $\Delta AEB$ and $\Delta BED$ right triangles $|AB|^2 - |EA|^2 = |BD|^2 - |ED|^2 \Rightarrow$

$$\Rightarrow a^2 - |EA|^2 = (a+1)^2 - (|EA|+2)^2 \Leftrightarrow |EA| = \frac{2a-3}{4}.$$

Similarly $|DF| = \frac{2b-3}{4}$. $|EF| = |EA| + |AD| + |DF| = \frac{2a-3}{4} + 2 + \frac{2b-3}{4} = \frac{a+b+1}{2}$.

i.e. $|EF| = |BK| = |BC|$. We will obtain that *ABCD* is rectangle. Thus

$$|BE| = |CF| \Rightarrow \frac{\sqrt{12a^2 + 12a - 9}}{4} = \frac{\sqrt{12b^2 + 12b - 9}}{4} \Leftrightarrow a(a+1) = b(b+1) \Leftrightarrow a = b.$$

I.e. $c = \frac{a+a+1}{2} = \frac{2a+1}{2} \notin N$.

This is impossible. Thus such Diophantine rectangle is not exist.

**Task 4**.

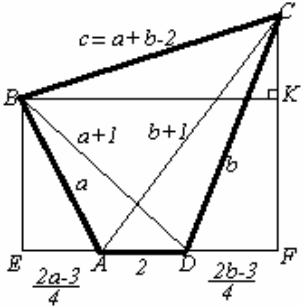

Is given: □ *ABCD* is circumscribed in the circle;

$$\begin{cases} |AB| = a; |BD| = a+1; |CD| = b; \\ |AC| = b+1; |AD| = 2; |BC| = c; \\ a, b, c \in N; b > a. \end{cases}$$

(In □ *ABCD* is inscribed the circle).

$$\Rightarrow c + 2 = a + b \Rightarrow c = a + b - 2$$

$$|EF| = \frac{2a-3}{4} + 2 + \frac{2b-3}{4} = \frac{2a+2b+2}{4} = \frac{a+b+1}{2}$$

$$\sqrt{c^2 - \frac{(a+b+1)^2}{4}} = |CF| - |BE| = \frac{\sqrt{12b^2 + 12b - 9}}{4} - \frac{\sqrt{12a^2 + 12a - 9}}{4} \Leftrightarrow$$

$$\Leftrightarrow 2\sqrt{4c^2 - (a+b+1)^2} = \sqrt{12b^2 + 12b - 9} - \sqrt{12a^2 + 12a - 9} \Leftrightarrow$$

$$\Leftrightarrow 4\left(4(a+b-2)^2 - (a+b+1)^2\right) = 12b^2 + 12b - 9 + 12a^2 + 12a - 9 - $$
$$- 2\sqrt{(12b^2 + 12b - 9)(12a^2 + 12a - 9)} \Leftrightarrow$$

$$\Leftrightarrow \sqrt{(12b^2 + 12b - 9)(12a^2 + 12a - 9)} = -39 - 12ab + 42a + 42b \Leftrightarrow$$

$$\Leftrightarrow \sqrt{(4b^2+12b-3)(4a^2+4a-3)} = -13-4ab+14a+14b \Leftrightarrow$$

$$\Leftrightarrow (4b^2+12b-3)(4a^2+4a-3) = (-13-4ab+14a+14b)^2 \Leftrightarrow$$

$$\Leftrightarrow 16a^2b^2+16ab^2-12b^2+16a^2b+16ab-12b-12a^2-12a+9 =$$
$$= 169+16a^2b^2+196a^2+196b^2+104ab-364a-364b-112a^2b-112ab^2+392ba \Leftrightarrow$$

$$\Leftrightarrow 128ab^2-208b^2+128a^2b-480ab+352b-208a^2+352a-160=0 \Leftrightarrow$$

$$\Leftrightarrow 8ab^2-13b^2+8a^2b-30ab+22b-13a^2+22a-10=0 \Leftrightarrow$$

$$\Leftrightarrow 8ab^2+8a^2b+22b+22a = 13b^2+13a^2+30ab+10 \quad (*)$$

If $\begin{cases} b \geq a > 3; \\ a,b \in N. \end{cases}$ then $\begin{cases} 8ab^2 \geq 32b^2 > 26b^2 \geq 13a^2+13b^2; \\ 8a^2b = 8a(ab) \geq 32ab > 30ab; \\ 22a+22b > 22 \cdot 3+22 \cdot 3 > 10. \end{cases} \Rightarrow$

$$\Rightarrow 8ab^2+8a^2b+22a+22b > 13a^2+13b^2+30a+10.$$

I.e. $\begin{cases} a,b \in N; \\ b \geq a; \\ a \in \{1;2;3\}. \end{cases}$

Due the introduction of $a = 1$, $a = 2$ and $a = 3$ in (*) we will obtain:

$$\begin{cases} \begin{bmatrix} 8b^2+8b+22b+22 = 13+13b^2+30b+10; \\ 16b^2+32b+22b+44 = 52+13b^2+60b+10; \\ 24b^2+72b+22b+66 = 117+13b^2+90b+10. \end{bmatrix} \\ b \in N. \end{cases} \Leftrightarrow \begin{cases} \begin{bmatrix} 5b^2+1=0; \\ b^2-2b-6=0; \\ 13b^2+4b-61=0. \end{bmatrix} \\ b \in N. \end{cases} \Leftrightarrow (b \in \varnothing).$$

I.e. such Diophantine tetragon is not exist.

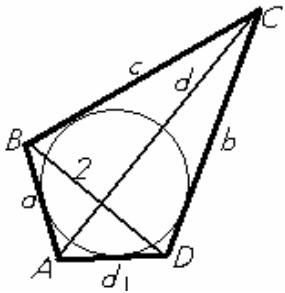

**Task 5**.  Is given: $ABCD$ is circumscribed on the circle;

$$\begin{cases} |CD|=b; |AD|=d_1; |AB|=a; |BC|=c; |AC|=d; \\ a+b=c+d; \ (1) \\ \begin{bmatrix} |BD|=2; \\ d_1=2. \end{bmatrix} \\ a,b,c \in N \end{cases}.$$

1) If $|BD| = 2$, then due the inequality of triangles and consideration of (1) we have:

**I.** $\begin{cases} b = c - 1; \\ d_1 = a - 1. \end{cases}$ ; **II.** $\begin{cases} b = c; \\ d_1 = a. \end{cases}$ ; **III.** $\begin{cases} b = c + 1; \\ d_1 = a + 1. \end{cases}$

Let's consider each of them

**I.**

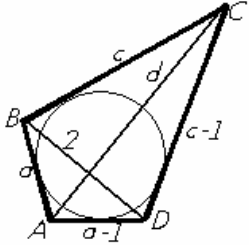

(From $\square ABCD$)  (From $\triangle ADC$)  $\begin{cases} 2 + d > c + a - 1; \\ d < (c-1) + (a-1); \\ a, c, d \in N. \end{cases} \Rightarrow \begin{cases} c + a - 3 < d < c + a - 2 \\ a, c, d \in N. \end{cases}$

This is impossible.

**II.**

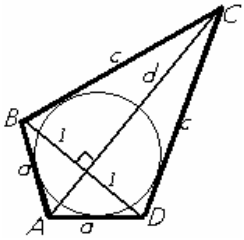

This case is considered in I of CHAPTER III.

**III.**

(From $\square ABCD$)  (From $\triangle ABC$)  $\begin{cases} d + 2 > (a+1) + c; \\ d < a + c; \\ a, c, d \in N. \end{cases} \Rightarrow \begin{cases} a + c - 1 < d < a + c; \\ a, c, d \in N. \end{cases}$

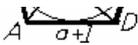

This is impossible.

I.e. does not existing the circumscribed on the circle Diophantine tetragon those length of diagonal is equal to 2.

**Remark**: A case when the length the diagonal of inscribed in the circle Diophantine tetragon is equal to 2 is considered in Task 7.

As already was mentioned, mathematicians were interested by the study of integer geometric figures since BC, despite this only some of the tasks occurs in this scope as oasis in the Gobi desert. The reason for this result is the absence of a unified method that would give the possibility us to study the properties of each Diophantine geometrical equation from the characteristic equations of this figure in a refined laconic way, to obtain in as simple as possible, the Diophantine equation (system of equations) andsimultaneously description in the optimal way for solving the equation. In this paper we will consideralso this issue.

We have proved in § 3 that there is no Diophantine n-gon (n> 3), both convex and concave, those lengths of arbitrary side or diagonal are equal to 1

Therefore, **the task * (n; 1)** is completely solved.

As we have seen above, solving the **task *(n; 2)** for k = 2 is the quiet complex problem even for n = 4, or for case of rectangles. In particular, as we will see below, even in the simplified case when the two sides of the Diophantine tetragon are parallel (the other two are parallel or not), the **task *(n; k)** for k = 2 is rather complex.

**Task 6. Find all Diophantine parallelograms and Diophantine trapezoids those length of one side or diagonal is equal to 2.**

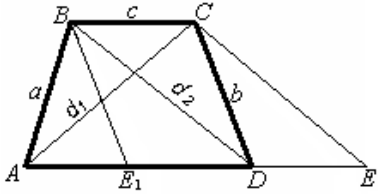

Fig. 6

It is known that for each Diophantine tetragon the length of each side and diagonal is more than 1, and if in Diophantine $\triangle MNP$ $|MN|=2$ and $|MP|=|NP|+k$, then $k \in \{-1; 0; 1\}$.

Let's say that conditions of task is satisfied by $\triangle ABCD$: $[BC]\|[AD]$; $|AB|=a$; $|CD|=b$; $|BC|=c$; $|AD|=d$; $|AC|=d_1$; $|BD|=d_2$. $a,b,c,d,d_1,d_2 \in N \setminus \{1\}$.

Let's drawn up: $[CE]\|[BD]$; $[BE_1]\|[CD]$; $E_1 \in [AD]$; $D \in [AE]$.

Without the limiting of generality let's say $a \leq b$ and $d \geq c$. It is easy to show that $d_1^2 + d_2^2 = a^2 + b^2 + 2cd$ (1) and $d_2 \geq d_1$ (2)

Let's consider the following cases:

I. $c=2$. From $\triangle ABC$ and $\triangle BDC$ we have: $d_1 \in \{a-1; a; a+1\}$; $d_2 \in \{b-1; b; b+1\}$.

In this case for solution of issue will be considered 9 cases, and as we will see further for $a = 2$ we will consider 27 cases. This is so labor consumptive that even at the solution of most complex problem, will be lost any desire to continue the work. There we shall propose the method that drastically reduces the number of considered at solution of this type tasks variants.

Let's say $d_1 = a+k$ and $d_2 = b+t$ (3). It is obvious that $k \in \{-1; 0; 1\}$ and $t \in \{-1; 0; 1\}$. Due the introduction of (3) in (1) $(a+k)^2 + (b+t)^2 = a^2 + b^2 + 4d$, from that by simplification we will obtain

$$2(2d - ak - bt) = b^2 + t^2.\qquad(4)$$

It is obvious that $k^2 + t^2 \neq 0$. Therefore from (4) due the taking into account that $k,t \in \{-1; 0; 1\}$, we have: $k^2 = t^2 = 1$. i.e. $k,t \in \{-1; 1\}$. Due the introduction the (4) we will obtain the four equations:

$$2d+a+b=1\,(5);\ 2d-a+b=1\,(6);\ 2d+a-b=1\,(7);\ 2d-a-b=1\,(8).$$

From $\triangle ABE_1$: $(d-c)+a \geq b$; $(d-c)+b \geq a$ (equalities take place when $d=c$ или $a=b$). I.e. $d+a-b \geq 2$ and $d-a+b \geq 2$, or from (5)-ь, (6) and (7), by given linmitations, has not the solution in natural numbers set. As for (8), in that case $k=t=1$. I.e. $d_1=a+1$ and $d_2=b+1$. It is easy to show that

$$(d+c)(a-b)(a+b) = (d-c)(d_1-d_2)(d_1+d_2).\qquad(9)$$

(This formula is valid as for $ABCD$ trapezoid, as well as for $ABCD$ parallelogram). If we consider that: $c=2$; $d_1=a+1$; $d_2=b+1$ and $a+b=2d$, then from (9) we will obtain:

$$(d+2)(a-b)(2d-1)=(d-2)(a-b)(2d+1) \Leftrightarrow$$

$$\Leftrightarrow (a-b)\bigl((d+2)(2d-1)-(d-2)(2d+1)\bigr)=0 \Leftrightarrow 6(a-b)d=0 \Leftrightarrow a=b \text{, but in that case they left part}$$

of (8) will be even, and right part the odd that is impossible. I.e. $c>2$.

II. From $2=a<c\le d$ $\triangle ABC$, $\triangle ABE_1$ and $\triangle ABD$ we will obtain that:

$$d_1 \in \{c-1;\ c;\ c+1\};\ b \in \{d-c-1;\ d-c;\ d-c+1\};\ d_2 \in \{d-1;\ d;\ d+1\}. \tag{10}$$

$\angle ABC$ and $\angle BAD$ are not acute, thus $d_1=c+1$ or $d_2=d+1$. If we consider that from $\triangle BCD$ $b+c>d_2$ and $b+c \in \{d-1;d;d+1\}$, we will obtain that $d_2<d+1$. I.e. we have $d_1=c+1$. I.e. (10) will be as following:

$$2=a<c\le d;\ d_1=c+1;\ b\in\{d-c-1;\ d-c;\ d-c+1\};\ d_2\in\{d-1;\ d\}. \tag{11}$$

Similarly to I let's say that $b=d-c+m$ and $d_1=d+n$ (12). It is obviously that: $m\in\{-1;\ 0;\ 1\}$ and $n\in\{-1;\ 0\}$ (13). Due the introduction of (12) in (1) we will obtain:

$$(c+1)^2+(d+n)^2=2^2+(d-c+m)^2+2dc \Leftrightarrow 2(dm-cm-c-dn+2)=n^2-m^2+1. \tag{14}$$

i.e. $n^2-m^2+1\equiv 0 \pmod 2$ (15). From (13) and (15) we have: $n=0$; $m\in\{-1;1\}$ (16) (if $n^2-m^2+1=2$, then $d=c-1$, that is impossible. I.e. $n^2-m^2+1=0$, from that will be obtained (16)). I.e.

$$\begin{cases} n^2-m^2+1=0; \\ 2(dm-cm-c-dn+2)=0. \end{cases} \Rightarrow \begin{cases} n=0, m\in\{-1;\ 1\}; \\ dm-cm-c+2=0. \end{cases} \Rightarrow \begin{bmatrix} -d+c-c+2=0; \\ d-c-c+2=0. \end{bmatrix} \Rightarrow \begin{bmatrix} d=2; \\ d=2c-2. \end{bmatrix}$$

$d=2$ is impossible, because $d>c>a=2$.

For $d=2c-2$ we will obtain: $b=c-1$; $d_2=2c-2$; $d_1=c+1$; $a=2$. Due introduction of these data in (9) we will obtain:

$$\begin{cases} (3c-2)(3-c)(c+1)=(c-2)(3-c)(3c-1); \\ c\in N\setminus\{1\}. \end{cases} \Leftrightarrow \begin{cases} \begin{bmatrix} c=\dfrac{1}{2}; \\ c=3. \end{bmatrix} \\ c\in N\setminus\{1\}. \end{cases} \Leftrightarrow c=3.$$

Thus we will obtain expressed on XII$_3$ trapezoid. (See Fig.7). (Let's notice that t6hhis was found (in 1989) the first Diophantine tetragon, those length of side is equal to 2).

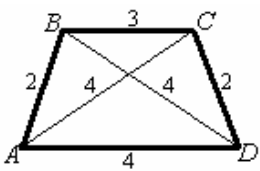

Fig. 7

III. $d_2>d_1=2$. Then from $\triangle ABC$, $\triangle BDC$ and $\triangle ACE$ we have: $a\in\{c-1;\ c;\ c+1\}$; $b\in\{d-1;\ d;\ d+1\}$; $d_2\in\{d+c-1;\ d+c\}$. In addition $d_2\ne d+c+1$ because $d_2<b+c\le d+c+1$. There also as in II let's say: $a=c+l$; $b=d+p$ and $d_2 = d + c + k$ (17). Obviously there also $l, p \in \{-1; 0; 1\}$; $k \in \{-1; 0\}$; (18). Due the introduction of (17) in (1) we will obtain:

$$2^2 + (d+c+k)^2 = (c+l)^2 + (d+p)^2 + 2cd \Leftrightarrow 2(ck + dk - lc - dp + 2) = l^2 + p^2 - k^2. \quad (19)$$

From (19) we have:

$$l^2 + p^2 - k^2 \equiv 0 (\mathrm{mod}\, 2) \Leftrightarrow \left( \{l = p = k = 0, \text{ an } \begin{cases} l = 0; \\ p = \pm 1; \\ k = -1. \end{cases} \text{ an } \begin{cases} p = 0; \\ l = \pm 1; \\ k = -1. \end{cases} \text{ an } \begin{cases} k = 0; \\ l = \pm 1; \\ p = \pm 1. \end{cases} \right). \quad (20)$$

For these values in the first three cases $ck + dk - lc - dp + 2 = 0$ (21$_1$), and in the fourth case $ck + dk - lc - dp + 2 = 1$ (21$_2$).

Due the try we will easy confirm that (21$_1$) and (21$_2$), accordingly of conditions (20), does not exist.

Finally we have that such Diophantine parallelogram does not exist, and from the Diophantine trapezoid only one satisfies to the conditions of task (see Fig.7). Let's consider once more task.

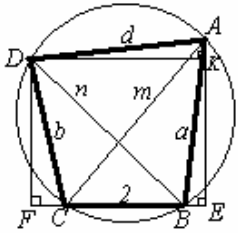
Fig. 8

**Task 7.** Let's find all such inscribed in circle Diophantine tetragon, those length of side or diagonal is equal to 2.

Let's say that $\triangle ABCD$ satisfies the condition of task and the length of it's any side without limitation of generality will be:

$$c = 2; \quad m = a + k; \quad n = b + l. \quad (22)$$

Obviously $k, l \in \{-1; 0; 1\}$. Accordingly of Ptolemy theorem (23) we have: $mn = ab + cd$ (24). Due the introduction of (22) in (24) and simplification we will obtain $al + bk + lk = 2d$ (25). Due the taking into account of (23) from (25) we have: $-b - a + 1 = 2d$ (26); $-b + a - 1 = 2d$ (27); $-b = 2d$ (28); $-a = 2d$ (29); $0 = 2d$ (30); $a = 2d$ (31); $b - a - 1 = 2d$ (32); $b = 2d$ (33); $b + a + 1 = 2d$ (34).

Obviously by stated in (26), (28) (29) and (30) limitations, have not the solutions in natural numbers. Also (27) has not the solution in natural numbers because in this case from $\triangle ABD$ $|AB| = a < |BD| + |AD| = = b + 1 + d < b + 1 + 2d$. Similarly we will obtain that (32) has not the solution in natural numbers. Let's now show that by stated limitations (34) also has not the solution in natural numbers. Let's assure the contrary. Let's say $\triangle ABCD$ satisfies the condition of task. In this case we have: $|CD| = b$; $|AD| = d = \dfrac{a+b+1}{2}$; $|AB| = a$; $|BC| = 2$; $|AC| = a+1$; $|DB| = b+1$. Let's plot: $E, F \in (CB); [DF] \perp (CB) \perp [AE]; K \in [EA) \perp [DK]$.

From the right triangles $\triangle AEB$ and $\triangle AEC$ $|AB|^2 - |BE|^2 = |AE|^2 = |AC|^2 - |CE|^2$. I.e. $a^2 - |BE|^2 =$

$= (a+1)^2 - (2+|BE|)^2 \Leftrightarrow |BE| = \dfrac{2a-3}{4}$. Similarly $|FC| = \dfrac{2b-3}{4}$. $|EF| = |FC| + |CB| + |BE| =$

$= \dfrac{2a-3}{4} + 2 + \dfrac{2b-3}{4} = \dfrac{a+b+1}{2} = |AD|$. I.e. $|DK| = |AD|$ თუ $|AK| = 0$ or $|DF| = |AE|$.

$\Rightarrow \sqrt{|AB|^2 - |BE|^2} = \sqrt{|CD|^2 - |FC|^2} \Leftrightarrow \sqrt{a^2 - \left(\dfrac{2a-3}{4}\right)^2} = \sqrt{b^2 - \left(\dfrac{2b-3}{4}\right)^2} \Rightarrow a(a+1) = b(b+1) \Rightarrow a = b$.

i.e. $d = \dfrac{a+b+1}{2} = \dfrac{2a+1}{2} \notin N$.

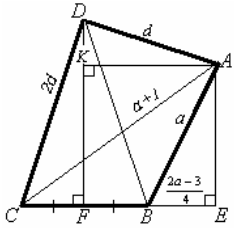

Fig. 9

This is impossible. I.e. our assumption is false or by stated limitations (34) has not the solution in natural numbers. Remains to be considered the identical equations (31) and (33).

$|AE| = 1 + \dfrac{2a-3}{4}$; $|DK| = \sqrt{d^2 - \left(\dfrac{2a+1}{4}\right)^2}$. $|AE| = |KF| = \sqrt{\dfrac{12a^2 + 12a - 9}{4}}$.

$\sqrt{|DB|^2 - |FB|^2} = |DF| = |DK| + |KF| = \sqrt{d^2 - \left(\dfrac{2a+1}{4}\right)^2} + \sqrt{\dfrac{12a^2 + 12a - 9}{4}}$

hence we have: $\left(\dfrac{3d^2}{a^2+a} - 1\right)^2 = 0 \Leftrightarrow \left((2a+1)^2 - 12d^2 = 1\right)$ (*).

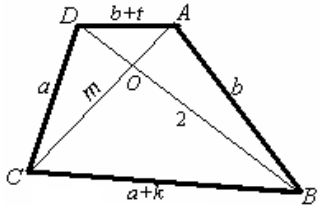

Fig. 10

This case reduces to considered by us XIV. Now let's say that the circle is inscribed in Diophantine $\triangle ABCD$ those lengths of diagonal without limitations let's say is:

$|DB| = 2$; $|AB| = b$; $|BC| = a+k$; $|CD| = a$; $|AD| = b+t$.

Is obvious that $t, k \in \{-1; 0; 1\}$ (35). Due the Ptolemy theorem we have: $2m = ab + (a+k)(b+t)$ (36) $\Leftrightarrow 2m = 2ab + at + bk + kt$.

From $\triangle ADC$ $m < a + b + t$

I.e. $2(a+b+t) > 2ab + at + bk + kt \Leftrightarrow 2ab + a(t-2) + b(k-2) + t(k-2) < 0$.

$\begin{cases} a > 2; \\ b > 2. \end{cases} \Rightarrow a \geq 3 \geq 1{,}5 + \dfrac{4{,}5}{2b-3} = \dfrac{3b}{2b-3} \Rightarrow 2ab - 3a - 3b \geq 0$

I.e. $2ab - 3a - 3b \geq 2ab + a(t-2) + b(k-2) + t(k-2) \Leftrightarrow t(2-k) \geq a(t+1) + b(k+1) \geq 0$. (37)

The (37) due the (35) condition has not the solution in natural numbers. I.e. $(a-2)(b-2) = 0$ (37$_1$).

From $\triangle AOD$ and $\triangle ACB$ due the application of inequality of triangles is easy to show that

$|AC| + |BD| > |AB| + |CD|$ i.e. $m + 2 > a + b \Rightarrow m > a + b - 2 \Rightarrow \begin{cases} m > a; \\ m > b. \end{cases}$ i.e.

$$m \in \{a+1;\ a+2;\ b+1;\ b+2\}. \tag{38}$$

With taking into account (38) from (36) we have:

$$2 = b(a-2) + (a+k)(b+t);  \quad (39) \quad 4 = b(a-2) + (a+k)(b+t); \tag{40}$$

$$2 = a(b-2) + (a+k)(b+t);  \quad (41) \quad 4 = a(b-2) + (a+k)(b+t). \tag{42}$$

It is easy to show that these equations have not the solution in natural numbersif:

$$\begin{cases} a = 2;\ b \in N \setminus \{1\}; \\ k, t \in \{-1;\ 0;\ 1\}. \end{cases} \text{ъ б } \begin{cases} b = 2;\ a \in N \setminus \{1\}; \\ k, t \in \{-1;\ 0;\ 1\}. \end{cases}$$

I.e.finally we have that length of each diagonal of inscribed in circle Diophantine tetragon is more that 2 and all having equal to 2 side Diophantine tetragon by the solution of Diophantine $(2a+1)^2 - 12d^2 = 1$ (43) equation

$$(|AB| = a;\ |BC| = 2;\ |CD| = 2d;\ |AD| = d;\ |DB| = 2d;\ |AC| = a+1). \tag{44}$$

**Let's mention that from the found Diophantine tetragons those length of any side is equal to 2, all will be inscribed in circle. Currently is not known this type of Diophantine tetragon on that would not circumscribed on circle. Even worse is situation related to the Diophantine tetragons those length of any diagonal is equal of 2. We have shown that if they exist, on them will not be circumscribed the circle, but what happens is other case is unknown to us, although as we have seen from the discussion of Task-2 (I-XVII), we have made attempts.**

The method shown by us in Task-2 and Task-3 does indeed significantly reduce the number of variants under consideration. However, the method of solving the obtained Diophantine equations is clearly expressed. All of this, for **task*(n;k)**, generally would be summarized as follows:

Consider n-gon vertices for all Diophantine triangles thhhose distance between any two vertices is equal to $k$, and taking into account that the modulus of difference of the other two sides represents an element of $\{1;\ 2;\cdots;k-1\}$ set, let us express the length of one of the sides by the length of other added $t$, where $t \in \{-(k-1);\ -(k-2);\ldots;\ 0;\ 1;\ 2;\ldots;\ (k-1)\}$;

By introducing the characteristic equation(s) of a given n-gon of this data, we obtain all possible Diophantine equations;

The excluding from the obtained Diophantine equations of equations which that due the stated limitations do not have solution in natural numbers, is mainly achieved by the additional limit6ations to the Diophantine equations obtained by using the triangle inequalities and also by the result of the **task*(n;k)** we have already obtained $\{((k;3));((k;4));\cdots((k;n-1))\}$.

By $((((k;\ n))$ is denoted condition that satisfies by the Diophantine n-gon sides and diagonals when the length of any side is equal to $k$).

Let's return to the **task*(n;k)** for $k = 2$.

Finally, we have that for the task*(n;2) convex $n \in \{3;4;5\}$ n-gon (though none of the Diophantine n-gons are found for n = 5 yet. In our opinion, such a pentagon does not exist) and for concave n-gons (here also for $n = 5$ and $n = 6$ we have the same picture as for convex for n = 5).

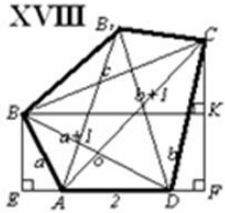
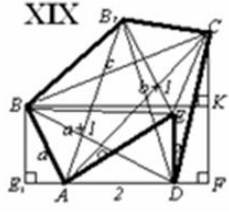
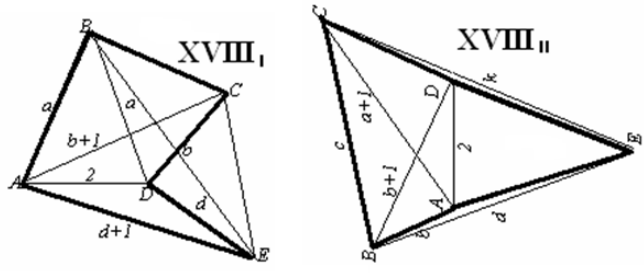

| Fig. 11 | Fig. 12 | Fig. 13 |

Let's now consider the task*(*n*;3).

Similarly to *k*=2, If the obe side of Diophantine triangle is equal to 3 then modulus of difference in other sides would be equal to 0, 1 or 2. With taking it into account all possible Diophantine tetragons for *k*-3 are presented below:

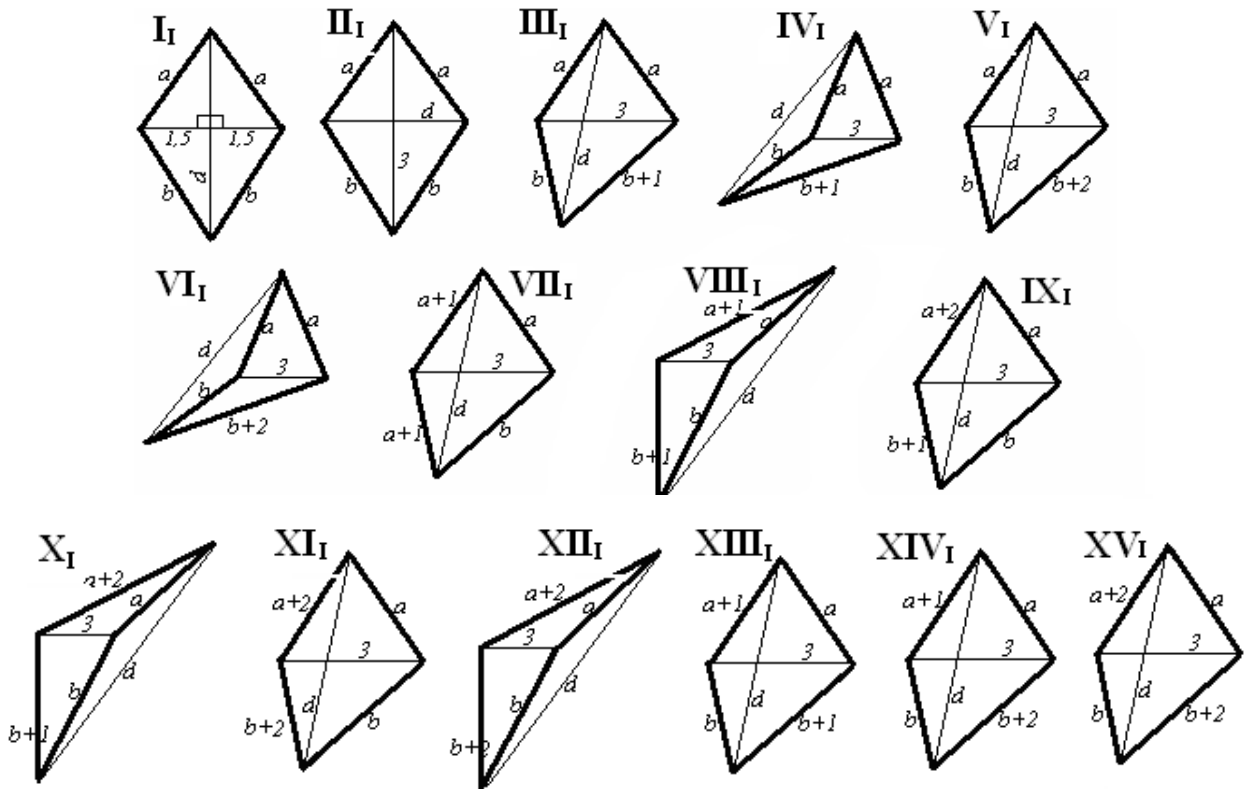

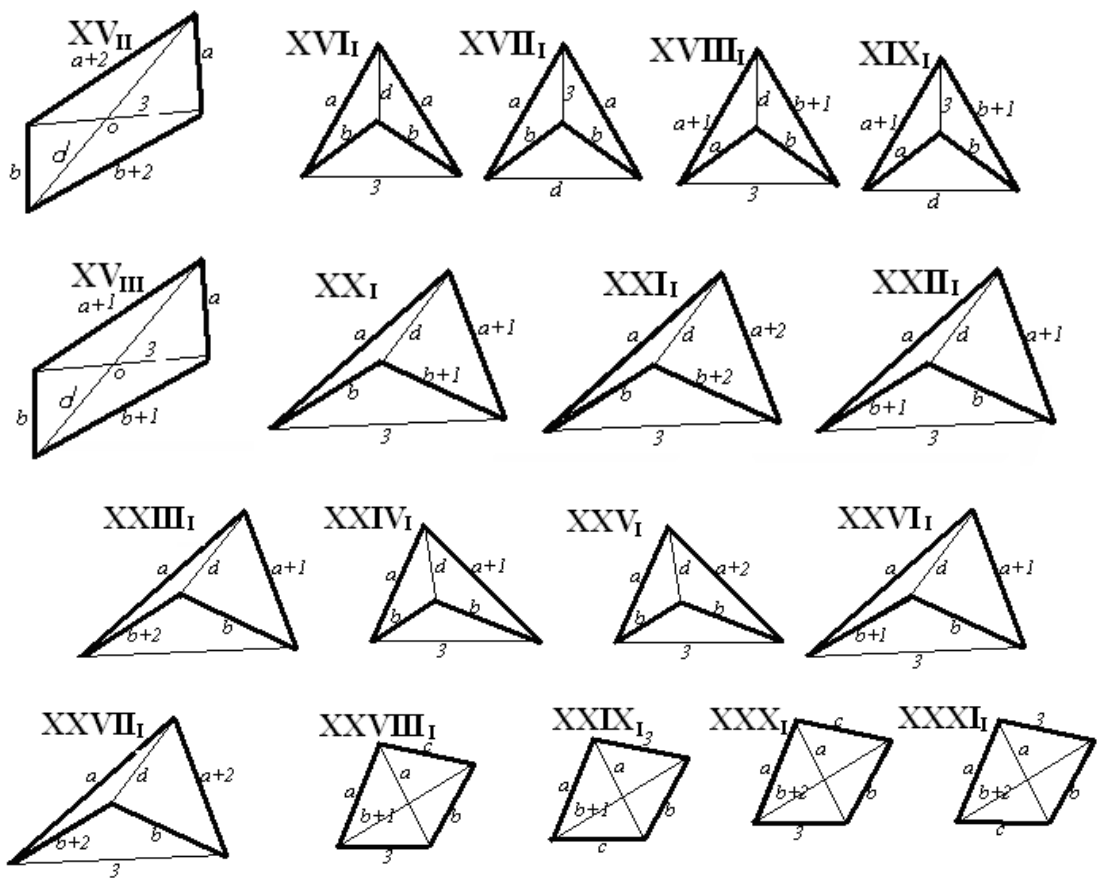

Let's consider from them the case when is obligatory for solution of issues of **task\***(*n*;3:

From the right triangles $\triangle AOB$ and $\triangle AOD$ we have:

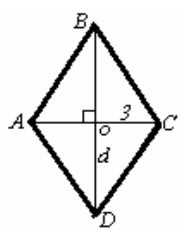

Is given: *ABCD*

$$\begin{cases} |AB|=|BC|=a; |AD|=|CD|=b; \\ |AC|=3; |BD|=d; \\ a,b,d \in N. \end{cases}$$

$$\begin{cases} 1.5^2 + (d-|OD|)^2 = a^2 \\ |OD|^2 + 1.5^2 = b^2 \\ a,b,d \in N \end{cases} \Leftrightarrow \begin{cases} |OD| = \dfrac{d^2+b^2-a^2}{2d} \\ \left(\dfrac{d^2+b^2-a^2}{d}\right)^2 + 9 = (2b)^2 \Rightarrow \\ \dfrac{d^2+b^2-a^2}{d} \equiv q \in N \end{cases}$$

$$\Rightarrow \begin{cases} (2b-q)(2b+q) = 9 \\ b,q \in N \end{cases} \Leftrightarrow \begin{cases} b = 2.5 \\ q = 4 \\ b,q \in N \end{cases} \Leftrightarrow (b \in \varnothing)$$

I.e. such Diophantine tetragon does not exist.

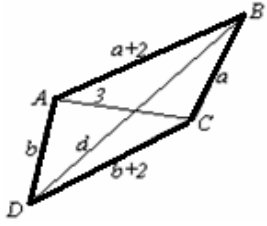

Is given: □ABCD

$$\begin{cases} |AD| = b; |AB| = a+2; |BC| = a; \\ |CD| = b+2; |AC| = 3; |BD| = d; \\ a,b,d \in N. \end{cases}$$

From $\triangle ABC$ and $\triangle ADC$ due the cosine law:

$$\begin{cases} a^2 = (b+2)^2 + 3^2 - 6(a+2) \cdot \cos B\widehat{A}C \\ (b+2)^2 = b^2 + 3^2 - 6b \cdot \cos C\widehat{A}D \end{cases} \Rightarrow \begin{cases} \cos B\widehat{A}C = \dfrac{2}{3} + \dfrac{1}{6(a+2)} \\ \cos C\widehat{A}D = -\dfrac{2}{3} + \dfrac{5}{6b} \end{cases} \Rightarrow$$

$$\Rightarrow \begin{cases} B\widehat{A}C \in \left(0; \arccos \dfrac{2}{3}\right) \\ C\widehat{A}D \in \left(0; \pi - \arccos \dfrac{2}{3}\right) \end{cases} \Rightarrow \left(B\widehat{A}D \in (0; \pi)\right)$$

(From $\triangle ABD$) $\quad \begin{cases} d < a+b+2 \\ d+3 > a+2+b+2 \\ a,b,d \in N \end{cases} \Rightarrow \begin{cases} a+b+1 < d < a+b+2 \\ a,b,d \in N \end{cases} \Rightarrow (d \in \varnothing).$

(From □ABCD)

I.e. such Diophantine tetragon does not exist.

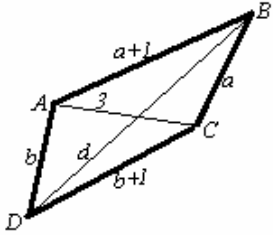

Is given: ABCD

$$\begin{cases} |AB| = a+1; |BC| = a; |CD| = b+1; \\ AD = b; |AC| = 3; |BD| = d; \\ a,b,d \in N. \end{cases}$$

There also as in previous task we have:

$$\begin{cases} \cos B\widehat{A}C = \dfrac{1}{3} + \dfrac{4}{3(a+1)} \\ \cos C\widehat{A}D = -\dfrac{1}{3} + \dfrac{4}{3b} \end{cases} \Rightarrow \begin{cases} B\widehat{A}C \in \left(0; \arccos \dfrac{1}{3}\right) \\ C\widehat{A}D \in \left(0; \pi - \arccos \dfrac{1}{3}\right) \end{cases} \Rightarrow C\widehat{A}D \in (0; \pi)$$

(From $\triangle ABD$) $\quad + b+1$

(From □ABCD) $\quad > a+1+b+1 \Rightarrow \begin{cases} a+b-1 < d < a+b+1 \\ a,b,d \in N \end{cases} \Rightarrow d = a+b.$

$\in N$

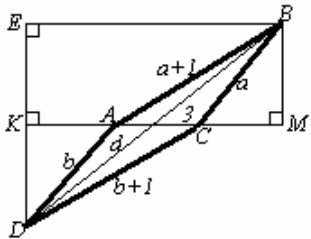

$$\begin{cases} b^2 - |KA|^2 = (b+1)^2 - (3+|KA|)^2 \\ b^2 - |CM|^2 = (a+1)^2 - (3+|CM|)^2 \end{cases} \Rightarrow$$

$$\Rightarrow \begin{cases} |KA| = \dfrac{b-4}{3}, \ |CM| = \dfrac{a-4}{3} \\ |DK| = \sqrt{b^2 - \left(\dfrac{b-4}{3}\right)^2} = \dfrac{\sqrt{8b^2+8b-16}}{3} \\ |BM| = |EK| = \dfrac{\sqrt{8a^2+8a-16}}{3} \end{cases}$$

From the $\triangle DEB$ right triangle $\left(|BD|^2 = |BE|^2 + |DE|^2\right)$, i.e

$$\begin{cases} (b+a)^2 = \left(\dfrac{a+b+1}{3}\right)^2 + \left(\dfrac{\sqrt{8a^2+8a-16} - \sqrt{8b^2+8b-16}}{3}\right)^2 \\ b, a \in N \end{cases} \Rightarrow$$

$$\Rightarrow \begin{cases} 2(8ab + 5a - 5b + 15) + 1 = 2\sqrt{(8a^2+8a-16)(8b^2+8b-16)} \\ a, b \in N \end{cases}.$$

This is impossible because the left side is even, and if right side is natural, then will be odd.

I.e. such Diophantine tetragon does not exist.

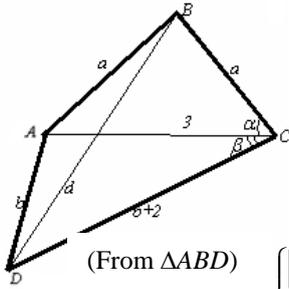

Is given: $ABCD$

$$\begin{cases} |AB| = |BC| = a; \ |CD| = b+2; \\ |AD| = b; |BD| = d; |AC| = 3; \\ a, b, d \in N. \end{cases}$$

(From $\triangle ABD$) $\begin{cases} |AC| + |BD| > |AB| + |DC|; \\ |AD| + |AB| > |BD|. \end{cases}$ $\Rightarrow \begin{cases} 3 + d > a + b + 2; \\ a + b > d; \\ a, b, d \in N \end{cases}$ $\Rightarrow \begin{cases} a + b - 1 < d < a + b; \\ a, b, d \in N \end{cases}$
(From $\square ABCD$)

This is impossible. i.e. such Diophantine tetragon does not exist.

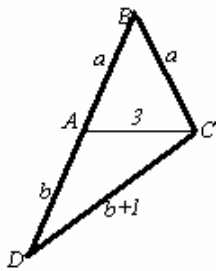

Is given: $\triangle DBC$

$$\begin{cases} A \in [BD]; |AB| = |BC| = a; \\ |AC| = 3; \\ |AD| = b; |CD| = b+1; \\ a, b, d \in N. \end{cases}$$

From $\triangle BAC$ and $\triangle DAC$ accordingly of cosine law we have:

$$\begin{cases} a^2 = a^2 + 3^2 - 6a\cos B\widehat{A}C; \\ (b+1)^2 = b^2 + 3^2 - 6b\cos(\pi - B\widehat{A}C); \Rightarrow \\ a, b \in N. \end{cases} \begin{cases} \cos B\widehat{A}C = \dfrac{3}{2a}; \\ \cos B\widehat{A}C = \dfrac{b-4}{3b}; \Rightarrow \\ a, b \in N. \end{cases}$$

$$\Rightarrow \begin{cases} \dfrac{3}{2a} = \dfrac{b-4}{3b}; \\ a,b \in N. \end{cases} \Rightarrow \begin{cases} b = 8 \cdot \dfrac{a}{2a-9}; \\ a,b \in N. \end{cases}$$

Because $(8; 2a-9) = 1$, thus $a \equiv 0 \,(\mathrm{mod}(2a-9)) \Rightarrow a \geq 2a - 9 \Rightarrow a \leq 9$.

We can find that only two pairs satisfies the conditions of task:

$$\begin{cases} a = 9 \\ b = 8 \end{cases} \text{and} \begin{cases} a = 5; \\ b = 16. \end{cases}$$

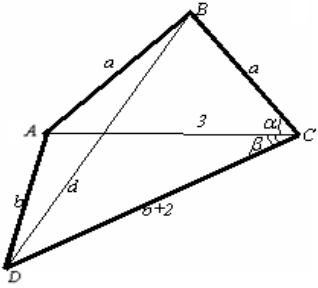

Is given: $\square ABCD$;

$$\begin{cases} |AB| = |BC| = a;\ |CD| = b+2; \\ |AD| = b;\ |BD| = d;\ |AC| = 3; \\ a,b,d \in N. \end{cases}$$

From $\triangle ABC$ and $\triangle DAC$ we have:

$$\begin{cases} \cos\alpha = \dfrac{3}{2a}; \\ (b+2)^2 = b^2 + 3^2 - 6b\cos\beta \end{cases} \Rightarrow \begin{cases} \cos\alpha \geq \dfrac{3}{4}; \\ \cos\beta = -\dfrac{2}{3} + \dfrac{5}{6b}. \end{cases} \Rightarrow \begin{cases} \alpha \in \left(\arccos\dfrac{3}{4}; \dfrac{\pi}{2}\right); \\ \beta \in \left(\dfrac{\pi}{2}; \dfrac{2\pi}{3}\right). \end{cases} \Rightarrow$$

$\Leftrightarrow \alpha + \beta \in \left(\dfrac{\pi}{2} + \arccos\dfrac{3}{4}; \dfrac{7\pi}{6}\right)$. i.e. $A$ and $C$ points related to $BD$ straight line would be located in one half-plane as well as in different half-planes.

The case when $d = a + b$

i.e. $\cos\widehat{BAC} = -\cos\widehat{CAD}$ or $\dfrac{3}{2a} = \dfrac{4b-5}{6b} \Leftrightarrow a = 2 + \dfrac{b+10}{4b-5}$.

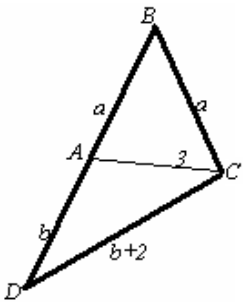

If $b > 5$, then $\dfrac{b+10}{4b-5} < 1$ and thus $a \notin N$.

If $b \leq 5$, then $\dfrac{b+10}{4b-5}$ is natural number when

$$\begin{cases} a = 3 \\ b = 2 \end{cases} \text{or} \begin{cases} a = 5 \\ b = 5 \end{cases}$$

From the above mentioned for the **Task*** $(n;3)\ 3 \leq n \leq 7$. In addition for $k=3$ is not found neither Diophantine pentagon nor Diophantine hexagon and Diophantine heptagon. In our opinion such Diophantine figures does not exis and in they exist their probable kind will be as followings:

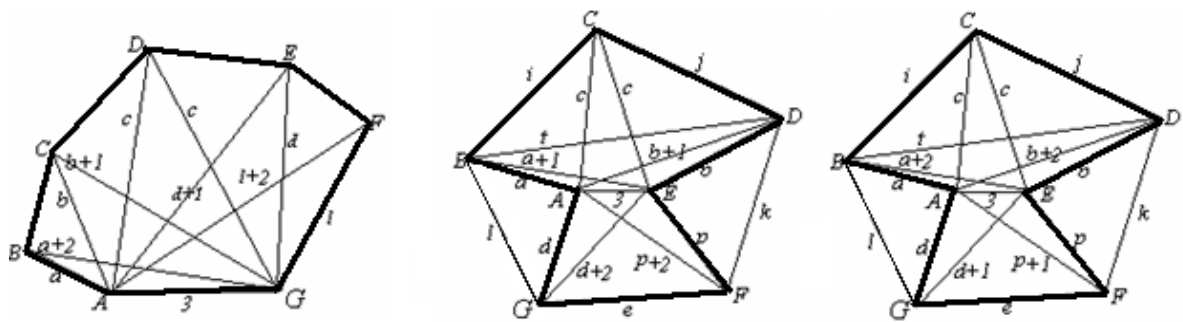

### ერთი ფუნდამენტური ამოცანის შესახებ დიოფანტურ გეომეტრიულფიგურებზე

ზურაბაღდგომელაშვილი

რეზიუმე

ნაშრომში დასმულია და შესწავლილია დიოფანტური $n$-კუთხედები. ამ სტატიის ავტორი დიოფანტურს უწოდებს მთელრიცხვა $n$-კუთხედს იმ მოტივით, რომ თითოეული მათგანის კომბინატორული თვისებების დასადგენად საჭიროა გარკვეული დიოფანტური განტოლების (განტოლებათა სისტემის) ამოხსნა. შემსწავლელი ამოცანებიდან ერთ-ერთი ფუნდამენტური ამოცანა.

**ამოცანა (n;k):** არსებობს თუ არა ყოველი ფიქსირებული $k$ნატურალური რიცხვისათვის ისეთი დიოფანტური $n$-კუთხედი ($n \geq 3$), რომლის რომელიმე გვერდის ან დიაგონალის სიგრძე ტოლია $k$-სი, და თუ არსებობს, მაშინ იპოვეთ ყველა ასეთი $n$.

ნაჩვენებია, რომ არ არსებობს ისეთი დიოფანტური $n$-კუთხედი ($n > 3$), როგორც ამოზნექილი, ასევე ჩაზნექილია, რომლის რომელიმე გვერდის ან დიაგონალის სიგრძე ტოლია ++-ის. i.e. $k = 1$-სათვის ზემოთხსენებული საკითხი გადაჭრილია.

კვლევის მეშვეობით ნაპოვნია რიგი დიოფანტური ოთხკუთხედები, რომელთა ერთ-ერთი გვერდი ტოლია 2-ის (აქ უნდა შევნიშნოთ, რომ ყოველი მათგანი წრეწირში ჩახაზული აღმოჩნდა). ნაჩვენებია, რომ ყოველი მათგანი წრეწირში ჩახაზული ისეთი დიოფანტური ოთხკუთხედი, რომლის დიაგონალის სიგრძე ტოლია k-ის. ნაჩვენებია, რომ ნებისმიერი ნატურალური ($k \geq 3$)-ისათვის მოიძებნება დიოფანტური ოთხკუთხედი, რომლის გვერდის სიგრძე ტოლია $k$-სი.

$n \in \{3;4;5\}$ ფუნდამენტური კვლევის შედეგად ნაჩვენებია, რომ $k = 2$-სათვის ამოზნექილი $n$-კუთხედებისათვის (თუმცა არც ერთი ასეთი დიოფანტური ხუთკუთხედი ჯერ არ არის ნაპოვნი. ავტორის აზრით ასეთი ხუთკუთხედი არ არსებობს) და $n \in \{3; 4; 5; 6\}$ – ჩაზნექილი $n$-კუთხედებისათვის (აქაც $n = 5$ და $n = 6$-სათვის არცერთი ასეთი $n$-კუთხედი არ არის ჯერ ნაპოვნი. და თუ არსებობს, შემთხვევისათვის მოცვანილია ასეთი ფიგურების ყველა შესაძლო სახე).

ნაშრომში განხილული ამოცანა (n;k), k=3-სათვის და ნაჭვენებია, რომ ამ შემთხვევაში $3 \le n \le 7$, მაგრამ არ არის ნაპოვნი არც ერთი ასეთი დიოფანტური ხუთკუთხედი, არც დიოფანტური ექვსკუთხედის და არც დიოფანტური შვიდკუთხედი, ავტორის აზრით ასეთი დიოფანტური ფიგურები არ არსებობენ, და თუ არსებობენ. მოცავს მათი სავარაუდო სახე.

# ОБ ОДНОЙ ФУНДАМЕНТАЛЬНОЙ ЗАДАЧЕ ПРОДИОФАНТОВЫМ ГЕОМЕТРИЧЕСКИМ ФИГУРАМ

## Зураб Агдгомелашвили

### Резюме

В труде поставлена и изучена одно из основным, ниже приведенных задач про диофантовым $n$-угольником (Автор статьи диофантовым называет целочисленным $n$-угольник, потому что, для установления комбинаторных свойства каждого из них требуется решить опредеденное диофантовон уравнение (систему уравнения).

**Задача (n;k):** Существует или нет для каждого фиксированного натурального числа такой диофантовый $n$-угольник $(n \ge 3)$, для которого длина некоторого стороны или диагонала равна $k$, и если существует, то найти все такие $n$.

Показана, что не существует такой диофантовый $n$-угольник $(n > 3)$, т.е. для $k = 1$ вышепоставленная задача решена.

Тщательным исследованием найдена ряд диофантовые четырехугольников сторона которых равна 2-м (Все они оказались вписанными в окружность). Показана, что не существует такой вписанный четырехугольник длина диагоналя которой равна 2-м. Показана что для любой $k \in N$ ($k \ge 3$) существует четырехугольник сторона которой равна $k$. Фундаментальным исследованием показана, что при $k = 2$, $n \in \{3; 4; 5\}$ длявыпуклых $n$-угольников (но пока не найдено ни одно такой диофантовый пятиугольник). По мнению автора такой пятиугольник не существует) и $n \in \{3; 4; 5; 6\}$, для вогнутых $n$-угольников (аналогично и здесь для $n = 5$ и $n = 6$).

# On a Fundamental Task of Diophantine Geometric Figures

## Zurab Agdgomelashvili

### Abstract

The goal of the work is to take on and study one of the fundamental tasks studying Diophantine n-gons (the author of the paper considers an integral n-gon is Diophantine as far as determination of combinatorial properties of each of them requires solution of a certain Diophantine equation (equation sets)).

Task*($n; k$): is there a Diophantine n-gon ($n \geq 3$) with any side or diagonal equal to *k* for each fixed natural number *k*. In case it exists then let us find each such *n*.

It is shown that there is no such Diophantine n-gon ($n > 3$), neither convex nor concave, the length of any side or diagonal of which is equal to 1. It means that for k=1 the above mentioned task is solved.

The studies made it possible to find certain Diophantine rectangles, one of the sides of which is equal to 2 (it is noteworthy that all of them appeared to be inscribed in a circumcircle). The studies showed that there is no such Diophantine rectangle inscribed ina circumcircle, the diagonal length of which is equal to 2. It is shown that for any natural ($k \geq 3$) there is a Diophantine rectangle with the side length equal to *k*.

The fundamental studies showed that for k=2 and n = {3; 4; 5}for convex n-gons (though such Diophantine pentagon has not been found yet. In the author's opinion such pentagon does not exist) and n= {3; 4; 5; 6} for concave n-gons (here as well for n=5 and n=6 no such n-gon has been found yet and for the case of its existence all probable types of such figures are presented).

The paper considers task *($n; k$) for *k=3* and shows that in this case 3≤ n ≤7. However, neither such Diophantine pentagon, nor Diophantine hexagon, nor Diophantine heptagon have been found yet. In the author's opinion such Diophantine figures do not exist and in case they do then he presents the probable types of them.